\documentclass[12pt, reqno]{amsart}
\usepackage{amssymb}
\usepackage{mathrsfs}
\usepackage{mathpazo}
\usepackage[colorlinks=true, breaklinks=true]{hyperref}
\usepackage{bm}
\usepackage[latin1]{inputenc}
\usepackage{graphicx}
\newtheorem{theorem}{Theorem}

\newtheorem{corollary}[theorem]{Corollary}

\newtheorem{definition}[theorem]{Definition}

\newtheorem{lemma}[theorem]{Lemma}

\newtheorem{proposition}[theorem]{Proposition}
\newtheorem{remark}[theorem]{Remark}

\newcommand{\bea}{\begin{eqnarray}}
\newcommand{\eq}{\end{eqnarray}}
\newcommand{\eea}{\end{eqnarray}}
\newcommand{\bqn}{\begin{eqnarray*}}
\newcommand{\beaa}{\begin{eqnarray*}}
\newcommand{\eqn}{\end{eqnarray*}}
\newcommand{\eeaa}{\end{eqnarray*}}
\newcommand{\bpr}{\begin{proposition}}
\newcommand{\epr}{\end{proposition}}

\newcommand{\cal}{\mathcal}

\numberwithin{equation}{section}
\numberwithin{theorem}{section}

\begin{document}
\title{Time-inhomogeneous Gaussian stochastic volatility models: Large deviations and super roughness}
\author{Archil Gulisashvili}\let\thefootnote\relax\footnotetext{Department of Mathematics, Ohio University, Athens OH 45701; e-mail: gulisash@ohio.edu}
\date{}

\begin{abstract}
We introduce time-inhomogeneous stochastic volatility models, in which the volatility is described by a nonnegative function of a Volterra type continuous Gaussian process that may have very rough sample paths. The main results obtained in the paper are sample path and small-noise large deviation principles for the log-price process in a time-inhomogeneous super rough Gaussian model under very mild restrictions. We use these results to study the asymptotic behavior of binary barrier options, exit time probability functions, and call options.
\end{abstract}

\maketitle

 \vspace{0.2in}

\vspace{0.2in}
\\
\section{Introduction}\label{S:I}
In this paper, we establish sample path and small-noise large deviation principles for general time-inhomogeneous Gaussian stochastic volatility models. We also study the asymptotic behavior of exit time probability functions, 
binary barrier options, and call options associated with such models. The main protagonists in the present paper are super rough Gaussian
models. In such a model, the volatility is described by a nonnegative time-dependent function of a Gaussian process that has extremely rough sample paths.

Let $(\Omega,{\cal F},\mathbb{P})$ be a probability space carrying two independent standard Brownian motions $W$ and $B$. We consider time-inhomogeneous stochastic volatility models, in which the asset price process $S_t$, $t\in[0,T]$, satisfies the following stochastic differential equation:
\begin{equation}
dS_t=S_tb(t,\widehat{B}_t)dt+S_t\sigma(t,\widehat{B}_t)(\bar{\rho}dW_t+\rho dB_t),\quad S_0=s_0> 0,\quad 0\le t\le T.
\label{E:mood}
\end{equation}
In (\ref{E:mood}), $s_0$ is the initial price, $T> 0$ is the time horizon, $\rho\in(-1,1)$ is the correlation coefficient, and  
$\bar{\rho}=\sqrt{1-\rho^2}$. The functions $b$ and $\sigma$ are continuous functions on $[0,T]\times\mathbb{R}$ satisfying a special condition (see Assumption C below). The equation in (\ref{E:mood}) is considered on the filtered probability space 
$(\Omega,\mathcal{F},\{\mathcal{F}_t\}_{0\le t\le T},\mathbb{P})$, where $\{\mathcal{F}_t\}_{0\le t\le T}$ is the augmentation of the filtration generated by the processes $W$ and $B$ (see \cite{KaS}, Definition 7.2 in Chapter 2).
We will also use the augmentation of the filtration generated by the process $B$, and denote it by
$\{\mathcal{F}^B_t\}_{0\le t\le T}$. The stochastic process $\widehat{B}$ appearing in (\ref{E:mood}) is a Volterra Gaussian process adapted to the filtration 
$\{\mathcal{F}^B_t\}_{0\le t\le T}$. The functions $b$ and $\sigma$ and the process $\widehat{B}$ appearing in (\ref{E:mood}) will be called the drift function, the volatility function, and the volatility process, respectively. The name ``Gaussian stochastic volatility model" was coined in \cite{GVZ1,GVZ2}, where a special model with $\rho=0$, $b(u)=0$, and $\sigma(u)=|u|$, $u\in\mathbb{R}$, was studied. 

In the last years, Gaussian fractional stochastic volatility models have become increasingly popular. The volatility in a fractional model is described by a nonnegative function of a fractional Volterra Gaussian process. Important examples of such processes are fractional Brownian motion and the Riemann-Liouville fractional Brownian motion. More information about these processes will be provided in Section \ref{S:log}. Short surveys of Gaussian fractional models can be found in \cite{GaS1,GS,G1}. 

The sample paths of fractional Brownian motion and the Riemann-Liouville fractional Brownian motion belong to certain H\"{o}lder spaces 
(see Remark \ref{R:ku}). In the present paper, we are mostly interested in continuous Volterra Gaussian processes such that their sample paths do not belong to any H\"{o}lder space. We call such processes super rough. In Section \ref{S:log}, we discuss  
Gaussian processes, for which the canonical metric function satisfies two-sided estimates with respect to a given modulus of continuity. Such processes were studied in \cite{MV} by Mocioalca and Viens (see Proposition 1 in \cite{MV}). This proposition is formulated in Section 
\ref{S:log} (see Theorem \ref{L:MV}). A Gaussian process appearing in Theorem \ref{L:MV} is super rough if the corresponding modulus of continuity grows near zero faster than any power function. The previous statement follows from Lemma \ref{L:mmm}. An interesting example of a super rough Gaussian process is logarithmic Brownian motion (logBm) introduced in \cite{MV} (see Section \ref{S:log}).
A Gaussian stochastic volatility model, in which the volatility process is super rough, will be called a super rough model. To the best of our knowledge, super rough Gaussian stochastic volatility models have never been studied before.

In Section \ref{S:log}, we define two special classes of Gaussian models (see Definitions \ref{D:gam} and \ref{D:SaS}). 
The volatility in a model from the first class is described by the Wick exponential of a constant multiple of one of the processes appearing in Theorem \ref{L:MV}, while the volatility in a model from the second class is the absolute value of such a process. If the conditions in Lemma \ref{L:mmm} are satisfied, then the models belonging to any of those classes are super rough. For instance, if we use the process logBm to build the volatility, we obtain a super rough model. The structure of a model from the first class is similar to that of the rough Bergomi model introduced in \cite{BFG}. However, the rough Bergomi model is not super rough (see Remark \ref{R:RB}). Applications of the rough Bergomi model in finance can be found in \cite{BHT,JMM,JPS,ZLCL}. The models from the second class introduced in Section \ref{S:log} resemble the Stein-Stein model (see \cite{SS}). The volatility in the Stein-Stein model is represented by the Ornstein-Uhlenbeck process. The Stein-Stein model is one of the classical stochastic volatility models. More information about this model can be found in \cite{G}. The paper \cite{LiS} discusses the importance of the Stein-Stein model in finance. The models from the second class are also related to the Gaussian models introduced in \cite{GVZ1,GVZ2}. 

Sample path and small-noise large deviation principles go back to the celebrated work of Varadhan (see \cite{V1}) and Freidlin and Ventsel' (see \cite{FW}). More information about large deviations can be found in \cite{DZ,DS,V2}. We refer the reader to \cite{BC,CMD,FGG,Ph,R} for applications of large deviation principles in financial mathematics. Throughout the paper, we use the abbreviation LDP to refer to a large deviation principle.

One of the main results in the present paper is Theorem \ref{T:27} formulated in Section \ref{S:ldr}. In this theorem, a sample path LDP is established for time-inhomogeneous Gaussian stochastic volatility models under very mild restrictions on the drift function, the volatility function, and the volatility process. More precisely, we assume that the drift function and the volatility function are locally $\omega$-continuous (see Definition \ref{D:moc}), and the modulus of continuity associated with the volatility process satisfies Fernique's 
condition (see (\ref{E:f30})). Lemma \ref{L:mmm} established in Section \ref{S:log} shows that under the restrictions formulated above, Gaussian stochastic volatility models can be super rough. To the best of our knowledge, Theorem \ref{T:27} provides the first sample path LDP that can be applied to super rough models. Note that the results on sample path and small-noise LDPs for Volterra Gaussian models obtained in the earlier works \cite{CP,FZ,G1,G2,JPS} do not apply to super rough models. For example, the small-noise LDP that was established 
in \cite{FZ}, Theorem 4.5, deals with driftless time-homogeneous fractional Gaussian models, in which the volatility function is globally H\"{o}lder continuous and the volatility process is fractional Brownian motion. Such a model is not super rough since fractional Brownian motion is not a super rough process. A sample path LDP for the rough Bergomi model was established in \cite{JPS}. This model is not super rough. In \cite{G1}, Theorem 13, a small-noise LDP was established for the log-price process in a driftless Gaussian model, in which the volatility function is locally $\omega$-continuous, while the volatility process is a Volterra Gaussian process such that the $L^2$-modulus of continuity of its kernel function satisfies a H\"{o}lder condition. In \cite{G2}, a sample path LDP was obtained under the same restrictions as those used in \cite{G1} 
(see Theorem 2.9 in \cite{G2}), while in \cite{CP}, an LDP similar to that in \cite{G2}, Theorem 2.9, was established for a time-homogeneous Volterra Gaussian model with drift (see Sections 6.2 and 6.3 in \cite{CP}). However, in the LDP obtained in \cite{CP}, an unnecessary extra restriction is imposed on the drift function and the volatility function. It is assumed in \cite{CP} that the volatility function and the absolute value of the drift function grow near infinity not faster than a power function (see Assumption 6.4 in \cite{CP}). 
The models considered in \cite{CP,G1,G2} are not super rough. Indeed, for any of those models, the modulus of continuity associated with the volatility process is of power type. The previous statement contradicts the definition of a super rough process. We would also like to mention the papers \cite{G3,HJL,GGG}, in which sample path LDPs were established for models with reflection, randomized fractional volatility models, and certain non-Gaussian stochastic volatility models, respectively. In \cite{BFGHS}, an asymptotic
expansion of the rate function in a small-noise LDP and a corresponding call pricing formula were obtained for certain fractional Gaussian models (see Theorems 3.1 and 3.2  in \cite{BFGHS}).

We will next provide a brief overview of the structure of the present paper. In Section \ref{S:aux}, we introduce general time-inhomogeneous Gaussian stochastic volatility models (see Definition \ref{D:gsvm}). Section \ref{S:aux} also contains various definitions and auxiliary statements. Section \ref{S:log} deals with super rough Gaussian processes and super rough Gaussian stochastic volatility models.
In Section \ref{S:ldr}, we formulate the main result of the present paper (Theorem \ref{T:27}) and discuss the structure of its proof. We also derive several corollaries from this theorem (see Theorems \ref{T:2}, \ref{T:17}, and \ref{T:1}).
Section \ref{S:addik} is devoted to applications of the LDPs established in Section \ref{S:ldr}. In Theorem \ref{T:333}, an asymptotic formula is obtained for binary barrier options, while Theorem \ref{T:3333} 
provides a large deviation style formula for the exit time probability function associated with the log-price process in a Gaussian model.
Formulas similar to the formula in Theorem \ref{T:3333} appeared earlier in a different context in the works of Freidlin and Ventsel' 
(see \cite{FW,FV11,FV12}). In \cite{CP,G2}, asymptotic formulas for the exit time probability functions were obtained for less general Gaussian models than those considered in the present paper. In Section \ref{S:addik}, we also study small-noise asymptotic behavior of binary call options and call options
in out-of-the-money regime (see Theorems \ref{T:final} and \ref{T:tre}, respectively). 
In the last section of the paper (Section \ref{S:101}), we prove Theorem \ref{T:27} (see Subsections \ref{SS:rv} -- \ref{SS:conti}).
\section{Time-Inhomogeneous Gaussian Stochastic Volatility Models}\label{S:aux}
In this section, we gather various definitions which will be used throughout the paper. Our first goal is to introduce Volterra Gaussian volatility processes.
\begin{definition}\label{D:accept}
Let $K$ be a real function on $[0,T]^2$. The function $K$ is called an admissible kernel if the following conditions hold: 
(a)\,$K$ is Borel measurable on $[0,T]^2$;
(b)\,$K$ is Lebesgue square-integrable over $[0,T]^2$;
(c)\,For every $t\in(0,T]$, the slice function $s\mapsto K(t,s)$, $s\in[0,T]$, 
belongs to the space $L^2[0,T]$;
(d)\,For every $t\in(0,T]$, the slice function is not almost everywhere zero.
\end{definition} 
An admissible kernel $(t,s)\mapsto K(t,s)$  is called a Volterra type kernel if
$K(t,s)=0$, for all $s> t$. 
It will be assumed throughout the paper that the volatility process $\widehat{B}$ 
in the Gaussian model introduced in (\ref{E:mood}) can be represented as follows: 
\begin{equation}
\widehat{B}_t=\int_0^tK(t,s)dB_s,\quad 0\le t\le T,
\label{E:f11}
\end{equation}
where $K$ is an admissible Volterra type kernel on $[0,T]^2$. The restriction on the slice functions in Definition \ref{D:accept} implies that the random variable $\widehat{B}_t$ is well-defined for every $t\in[0,T]$. If the kernel $K$ is just square-integrable over $[0,T]^2$, then only almost all slice functions are in $L^2[0,T]$. It is clear that the process $\widehat{B}$ is adapted to the filtration 
$\{\mathcal{F}^B_t\}_{0\le t\le T}$. An additional restriction on the kernel $K$ will be imposed below
(see Assumption A). Under Assumption A, the process $\widehat{B}_t$, $t\in[0,T]$, is a continuous Gaussian process. 
The covariance function of this process is given by
$
C(t,s)=\int_0^TK(t,u)K(s,u)du,
$
for all $t,s\in[0,T]$. It follows from condition $(d)$ in Definition \ref{D:accept} that the process $\widehat{B}$ is nondegenerated in the following sense: The variance function of $\widehat{B}$ given by $V(t)=\int_0^TK(t,u)^2du$ satisfies $V(t)> 0$, for all $t\in(0,T]$.

Popular examples of Volterra Gaussian volatility processes are fractional Brownian motion and the Riemann-Liouville fractional Brownian motion. 
For $0< H< 1$, fractional Brownian motion $B^H_t$, $t\ge 0$, is a centered Gaussian process with the covariance function given by
$C_H(t,s)=\frac{1}{2}\left(t^{2H}+s^{2H}-|t-s|^{2H}\right)$, for all $t,s\ge 0$.
The process $B^H$ was first implicitly considered by Kolmogorov in \cite{Ko}, and was studied by Mandelbrot and van Ness in \cite{MvN}. The constant $H$ is called the Hurst parameter. If $H=\frac{1}{2}$, then we get standard Brownian motion. It is clear that the variance function of the process $B^H$ is given by $V^{H}(t)=t^{2H}$, $t\ge 0$. The Riemann-Liouville fractional Brownian motion
is defined as follows: $R^H_t=\frac{1}{\Gamma(H+\frac{1}{2})}\int_0^t(t-s)^{H-\frac{1}{2}}dB_s$, $t\ge 0$, 
where $H> 0$, and $\Gamma$ is the gamma function. 
The process $R^H$ is Gaussian, and its covariance function is given by 
$$
\widetilde{C}_H(t,s)=\frac{t^{H-\frac{1}{2}}s^{H+\frac{1}{2}}}
{(H+\frac{1}{2})\Gamma(H+\frac{1}{2})^2}\,_2F_1(\frac{1}{2}-H,1,H+\frac{3}{2},\frac{s}{t}),\quad 0< s< t,
$$
where $_2F_1$ is the hypergeometric function (see (11) in \cite{LS}). The variance function of $R^H$ is as follows:
$
\widetilde{V}(t)=\frac{t^{2H}}{2H\Gamma(H+\frac{1}{2})^2}, 
$
$t> 0$.
More information about the process $R^H$ can be found in \cite{LS,Pi}. The processes $B^H$ and $R^H$ are Volterra Gaussian processes. This is clear for the process $R^H$. For fractional Brownian motion, a Volterra type representation was found in \cite{MG} by Mol$\check{\rm c}$an and Golosov (see also \cite{DU}). The kernel $K_H$ in the Volterra type representation of the process $B^H$ is as follows. 
For $\frac{1}{2}< H< 1$, 
$$
K_H(t,s)=\sqrt{\frac{H(2H-1)}{\int_0^1(1-x)^{1-2H}x^{H-\frac{3}{2}}dx}}s^{\frac{1}{2}-H}
\int_s^t(u-s)^{H-\frac{3}{2}}u^{H-\frac{1}{2}}du,\quad 0< s< t,
$$
while for $0< H<\frac{1}{2}$,
\begin{align*}
&K_H(t,s)=\sqrt{\frac{2H}{(1-2H)\int_0^1(1-x)^{-2H}x^{H-\frac{1}{2}}dx}} \\
&\quad\left[\left(\frac{t}{s}\right)^{H-\frac{1}{2}}(t-s)^{H-\frac{1}{2}}-\left(H-\frac{1}{2}\right)
s^{\frac{1}{2}-H}\int_s^tu^{H-\frac{3}{2}}(u-s)^{H-\frac{1}{2}}du\right],\quad 0< s< t.
\end{align*}

Moduli of continuity will be broadly used throughout the present paper.
\begin{definition}\label{D:modc}
(i)\,Let $r> 0$. A bounded function $\eta:[0,r]\mapsto[0,\infty)$ is called a modulus of continuity on $[0,r]$,
if $\eta(0)=0$ and $\displaystyle{\lim_{u\rightarrow 0}\eta(u)=0}$;
(ii)\,A locally bounded function $\omega:[0,\infty)\mapsto[0,\infty)$ is called a modulus of continuity
on $[0,\infty)$, if $\omega(0)=0$ and $\displaystyle{\lim_{u\rightarrow 0}\omega(u)=0}$.
\end{definition}

Let $X_t$, $t\in[0,T]$, be a Gaussian process on $(\Omega,{\cal F},\mathbb{P})$. The canonical metric function $\delta$ associated with this process is defined by
\begin{equation}
\delta^2(t,s)=\mathbb{E}[(X_t-X_s)^2],\quad(t,s)\in[0,T]^2.
\label{E:cam}
\end{equation}
Suppose $\eta$ is a modulus of continuity on $[0,T]$ such that $\delta(t,s)\le\eta(|t-s|)$, $t,s\in[0,T]$,
and for some $b> 1$, the following inequality holds:
\begin{equation}
\int_{b}^{\infty}\eta\left(u^{-1}\right)(\log u)^{-\frac{1}{2}}\frac{du}{u}<\infty.
\label{E:f30}
\end{equation}
\begin{lemma}\label{L:ferni}
A Gaussian process $X_t$, $t\in[0,T]$, satisfying the previous two conditions is a continuous stochastic process. 
\end{lemma}

The statement in Lemma \ref{L:ferni} was announced by Fernique in \cite{Fer}. The first proof was published by Dudley in \cite{Dud}
(see also \cite{MS} and the references therein).

By the It$\hat{\rm o}$ isometry, the following equality holds 
for the process $\widehat{B}$:
\begin{equation}
\delta^2(t,s)=\int_0^T(K(t,u)-K(s,u))^2du,\quad t,s\in[0,T].
\label{E:iso1}
\end{equation}
The $L^2$-modulus of continuity of the kernel $K$ is defined on $[0,T]$ by
\begin{equation}
M_K(\tau)=\sup_{t,s\in[0,T]:|t-s|\le\tau}\int_0^T(K(t,u)-K(s,u))^2du,
\label{E:iso2}
\end{equation}
for all $\tau\in[0,T]$. If $\lim_{\tau\rightarrow 0}M_K(\tau)=0$, then the variance function of the process $\widehat{B}$ defined in 
(\ref{E:f11}), that is, the function $V(t)=\int_0^TK(t,s)^2ds$, $t\in[0,T]$,
is continuous, and hence
\begin{equation}
\max_{t\in[0,T]}\int_0^tK(t,s)^2ds<\infty.
\label{E:mmn}
\end{equation}
The continuity of the function $V$ can be established as follows. First, using the Cauchy-Schwarz inequality, we obtain the estimate 
$$
|V(t)^{\frac{1}{2}}-V(s)^{\frac{1}{2}}|\le\{\int_0^T(K(t,u)-K(s,u))^2du\}^{\frac{1}{2}},\quad t,s\in[0,T].
$$
Then, by taking into account the condition $\lim_{\tau\rightarrow 0}M_K(\tau)=0$ and the fact that $K$ is an admissible kernel, we finish the proof of the continuity of $V$.

We will next formulate an assumption that will be used in the rest of the paper.  \\
\\
\underline{Assumption A:} The $L^2$-modulus of continuity of the kernel $K$ satisfies the following condition:
$M_K(\tau)\le\eta(\tau)^2$, $\tau\in[0,T]$, where $\eta$ is a modulus of continuity on $[0,T]$ such that (\ref{E:f30}) holds.
\\
\\
Assumption A guarantees that the process $\widehat{B}$ in (\ref{E:f11}) is a continuous Gaussian process.
\begin{remark}\label{R:fty}
The kernel of fractional Brownian motion $B^H$ satisfies Assumption A with $\eta(\tau)=a_Ht^{2H}$, $a_H> 0$, $t\in[0,T]$.
The previous fact was established in \cite{Z} (see also Lemma 8 in \cite{G1}). Assumption A also holds for the kernel of the 
Riemann-Liouville fractional Brownian motion $R^H$ (see \cite{G1}, Lemma 8). Here $\eta(\tau)=b_Ht^{2H}$, $b_H> 0$, $t\in[0,T]$.
More examples of kernels satisfying the condition in Assumption A will be given in Section \ref{S:log}.
\end{remark}

Let $x=(u_1,v_1)$ and $y=(u_2,v_2)$ be elements of the metric space $[0,T]\times\mathbb{R}$ equipped with the distance
$d(x,y)=||x-y||=|u_1-u_2|+|v_1-v_2|$. The closed ball in this space having radius $\delta> 0$ and centered at $(0,0)$ 
will be denoted by $\overline{B(\delta)}$. We will next formulate certain restrictions on the drift function $(u,v)\mapsto b(u,v)$ and the volatility function 
$(u,v)\mapsto\sigma(u,v)$, where $(u,v)\in[0,T]\times\mathbb{R}$. 
\begin{definition}\label{D:moc}
Let $\omega$ be a modulus of continuity on $[0,\infty)$.
A function $\lambda$ defined on $[0,T]\times\mathbb{R}$ is called locally $\omega$-continuous,
if for every $\delta> 0$ there exists a number $L(\delta)> 0$ such that for all $x,y\in\overline{B(\delta)}$,
the following inequality holds:
$
|\lambda(x)-\lambda(y)|\le L(\delta)\omega(||x-y||).
$
\end{definition}
\begin{remark}\label{R:let}
Let us denote by $MC$ the class of all strictly increasing continuous moduli of continuity $\mu$ on $[0,\infty)$
such that $\lim_{\delta\rightarrow\infty}\mu(\delta)=\infty$.
With no loss of generality, we will assume throughout the rest of the paper that the modulus of continuity $\omega$ in Definition \ref{D:moc}
is from the class $MC$. Indeed, it is not hard to prove that any modulus of continuity $\omega$ can be majorized by a modulus of continuity from the class $MC$. In addition, it is not difficult to see that we can also assume that $L\in MC$, where $L$ is the function appearing in Definition \ref{D:moc}.
\end{remark}

The next condition restricts the class of admissible drift and volatility functions.
\vspace{0.1in}
\\
\underline{Assumption C:} The drift function $b$ and the volatility function $\sigma$ are locally $\omega$-continuous on the space 
$[0,T]\times\mathbb{R}$ for some modulus of continuity $\omega$ on $[0,\infty)$. The volatility function $\sigma$ is not identically zero
on $[0,T]\times\mathbb{R}$.
\\
\\
In Section \ref{S:log}, we provide two examples of volatility functions used in financial modeling of volatility. The first example is the volatility function employed in the rough Bergomi model (see Remark \ref{R:RB}), while the second example is the function $\sigma(u)=c|u|$, $c> 0$, $u\in\mathbb{R}$, used in the Stein-Stein model (see the discussion at the end of Section \ref{S:log}). It is not hard to see that these functions satisfy the conditions in Assumption C.
\begin{definition}\label{D:gsvm}
An asset price model described by the stochastic differential equation in (\ref{E:mood}) will be called a general time-inhomogeneous
Gaussian stochastic volatility model, if the following conditions hold: 
(a)\,The drift function $b$ and the volatility function $\sigma$ satisfy the conditions in Assumption C.\,\,
(b)\,Assumption A holds for the kernel $K$ of the volatility process $\widehat{B}$ defined by (\ref{E:f11}). 
\end{definition}
We have already mentioned that if condition (b) formulated above holds for the kernel $K$, then the volatility process 
$\widehat{B}$ is a continuous Gaussian process.  


Let $\varepsilon\in(0,1]$ be a small-noise parameter. We will work with the following scaled version of the model in (\ref{E:mood}):
\begin{equation}
dS^{(\varepsilon)}_t=S^{(\varepsilon)}_tb(t,\sqrt{\varepsilon}\widehat{B}_t)dt
+\sqrt{\varepsilon}S^{(\varepsilon)}_t\sigma\left(t,\sqrt{\varepsilon}\widehat{B}_t\right)(\bar{\rho}dW_t+\rho dB_t),
\label{E:bbh}
\end{equation}
where $0\le t\le T$. For every fixed $\varepsilon\in(0,1]$, the equation in (\ref{E:bbh}) is a linear stochastic differential equation with respect to the process 
\begin{equation}
t\mapsto\int_0^tb(s,\sqrt{\varepsilon}\widehat{B}_s)ds+\sqrt{\varepsilon}
\int_0^t\sigma\left(s,\sqrt{\varepsilon}\widehat{B}_s\right)(\bar{\rho}dW_s+\rho dB_s),\quad t\in[0,T].
\label{E:pre}
\end{equation} 
The process in (\ref{E:pre}) is a semimartingale. The previous statement follows from the inequalities  
$\int_0^T|b(s,\sqrt{\varepsilon}\widehat{B}_s)|ds<\infty$ and $\int_0^T\sigma\left(s,\sqrt{\varepsilon}\widehat{B}_s\right)^2ds<\infty$
$\mathbb{P}$-a.s. The unique solution to the equation in (\ref{E:bbh}) is the Dol\'{e}ans-Dade exponential
$$
S_t^{(\varepsilon)}=s_0\exp\left\{\int_0^tb(s,\sqrt{\varepsilon}\widehat{B}_s)ds
-\frac{\varepsilon}{2}\int_0^t\sigma^2(s,\sqrt{\varepsilon}\widehat{B}_s)ds
+\sqrt{\varepsilon}\int_0^t\sigma(s,\sqrt{\varepsilon}\widehat{B}_s)
(\bar{\rho}dW_s+\rho dB_s)\right\},
$$
where $0\le t\le T$ (see, e.g., \cite{RY}). We will denote by $X^{(\varepsilon)}$ the scaled log-price process given by $X_t^{(\varepsilon)}=\log S_t^{(\varepsilon)}$, $0\le t\le T$, 
and set $x_0=\log s_0$. 
\section{Super Rough Gaussian Stochastic Volatility Models}\label{S:log}
This section deals with Gaussian models, in which the volatility process has extremely rough sample paths. 
\begin{definition}\label{D:cf}
Let $\eta$ be a continuous increasing modulus of continuity on $[0,T]$. The function class $C^{(\eta)}$ is defined as the class of all continuous functions on $[0,T]$ that admit $\eta$ as a uniform modulus of continuity. More precisely, 
\begin{equation}
f\in C^{(\eta)}\Longleftrightarrow \sup_{0\le s< t\le T}\frac{f(t)-f(s)}{\eta(t-s)}<\infty.
\label{E:hol}
\end{equation}
If $r\in(0,1)$ and $\eta_r(t)=t^r$, $t\in[0,T]$, then the condition in (\ref{E:hol}) is the H\"{o}lder condition, and the corresponding 
function class will be denoted by $C^r$.
\end{definition}
\begin{remark}\label{R:ku}
It is known that $\mathbb{P}$-almost all sample paths of fractional Brownian motion $B^H$ belong to the class
$\bigcap_{0< r< H}C^{r}$, but do not belong to any of the classes $C^{H+\delta}$ with $\delta> 0$ (see, e.g., \cite{BHOZ}, Section 1.6).
The sample paths of the Riemann-Liouville fractional Brownian motion $R^H$, $0< H< 1$, have the same properties. Indeed, it follows from Theorem 5.1 in \cite{Pi} that the processes $R^H$ and $B^H$ can be jointly realized so that the process $R^H-B^H$ has infinitely differentiable paths.
\end{remark}
\begin{definition}\label{D:sr}
A Gaussian process for which $\mathbb{P}$-almost all sample paths do not belong to any of the classes $C^{r}$, $0< r< 1$, will be called a super rough Gaussian process. 
\end{definition}
It follows from the remarks made above that the processes $B^H$ and $R^H$ are not super rough.
The next assertion was established in \cite{MV} by Mocioalca and Viens. 
\begin{theorem}\label{L:MV}
Let $\eta$ be a modulus of continuity on $[0,T]$, and suppose
 $\eta\in\mathbb{C}^2(0,T)$. Suppose also that the function $x\mapsto(\eta^2)^{\prime}(x)$ 
is positive and non-increasing on $(0,T)$. Set
\begin{equation}
\tau(x)=\sqrt{(\eta^2)^{\prime}(x)},\quad x\in(0,T). 
\label{E:tr}
\end{equation}
Then the process
\begin{equation}
\widehat{B}_t^{(\eta)}=\int_0^t\tau(t-s)dB_s,\quad t\in[0,T],
\label{E:ty}
\end{equation}
is a Gaussian process satisfying the following conditions: (i)\,There exist positive constants 
$c_1$ and $c_2$ such that
$\eta(|t-s|)\le\delta(t,s)\le 2\eta(|t-s|)$, for all $t,s\in[0,T]$. Here the symbol $\delta$ stands for the canonical metric function associated with the process $\widehat{B}^{(\eta)}$; (ii)\,$\widehat{B}_0^{(\eta)}=0$\,$\mathbb{P}$-a.s.; 
(iii)\,The process $\widehat{B}^{(\eta)}$ is adapted to the filtration $\{\mathcal{F}^B_t\}_{0\le t\le T}$. 
\end{theorem}
\begin{remark}\label{R:oi}
A modulus of continuity $\eta$ such that the conditions in Theorem \ref{L:MV} hold for it does not necessarily satisfy Fernique's condition in Assumption A. Such an example is the logarithmic modulus of continuity $\eta_{\beta}$ with $0<\beta\le 1$ (see the definition in (\ref{E:logg})
below).
\end{remark}

Theorem \ref{L:MV} can be used to build super rough processes.
\begin{lemma}\label{L:mmm}
Suppose $\eta$ is a modulus of continuity satisfying the conditions in Theorem \ref{L:MV} and growing near zero faster than any power function $x\mapsto x^a$, $a> 0$. Suppose also that $\eta$ satisfies the conditions
in Assumption A. Then the process $\widehat{B}^{(\eta)}$ defined by (\ref{E:ty}) is a continuous super rough Volterra Gaussian process.  
\end{lemma}

\it Proof. \rm Since Assumption A holds, the process $\widehat{B}^{(\eta)}$ is continuous. We will next prove that this process is super rough. Using the restrictions on the modulus of continuity $\eta$ in Theorem \ref{L:MV}, we see that the function $\eta^2$ is 
concave on $(0,T)$. Moreover $\eta^2(0)=0$. It is known that such functions are subadditive. Therefore, for $0< s\le t< T$, we have
$\eta^2(t-s)\ge\eta^2(t)-\eta^2(s)$. Next, using (i) in Theorem \ref{L:MV}, we see that
\begin{equation}
\delta^2(t,s)\ge\eta^2(t)-\eta^2(s),\,\,0< s\le t< T.
\label{E:vid}
\end{equation}

Our next goal is to use Lemma 7.2.13 in \cite{MR} with $G=\widehat{B}^{(\eta)}$, $f=\eta$, and $c=1$. Condition (7.177) in 
Lemma 7.2.13 follows from (\ref{E:vid}). Applying the lemma and using the fact that $\eta$ is an increasing function, we obtain the following estimate:
$$
\limsup_{t\rightarrow 0}\frac{|\widehat{B}^{(\eta)}_t|}{\eta(t)\sqrt{\log\log\frac{1}{\eta^2(t)}}}\ge\sqrt{2}
$$
$\mathbb{P}$-a.s. Now, it is not hard to see that the previous estimate and the condition that $\eta$ grows near zero faster than any power function imply the super roughness of the process $\widehat{B}^{(\eta)}$.

This completes the proof of Lemma \ref{L:mmm}.

General time-inhomogeneous Gaussian stochastic volatility models were introduced in Definition \ref{D:gsvm}. In the remaining part of the present section, we will discuss special classes of such models. 
\begin{definition}\label{D:Dt}
A super rough Gaussian stochastic volatility model is a model,
in which the volatility process $\widehat{B}$ is a super rough continuous Volterra Gaussian process. 
\end{definition}
\begin{remark}\label{R:nbn}
Suppose the volatility process in a Gaussian model is the process $\widehat{B}^{(\eta)}$, where $\eta$ satisfies the conditions in 
Lemma \ref{L:mmm}. Then the model is super rough. The previous statement follows from Lemma \ref{L:mmm}.
\end{remark}

For a zero-mean Gaussian random variable $G$, the Wick exponential is defined by 
$
{\cal E}(G)=\exp\{-\frac{1}{2}\mathbb{E}[G^2]+G\}
$
(see, e.g., \cite{BP}, p. 392). For the square integrable Volterra-type kernel $K(t,s)=\tau(t-s)\mathbb{1}_{s\le t}$, $(t,s)\in[0,T]^2$, appearing in formula (\ref{E:ty}), the variance function of the corresponding process $\widehat{B}_t^{(\eta)}$ is $t\mapsto\eta(t)^2$, $t\in[0,T]$. Hence, 
for every $c> 0$ and $t\in[0,T]$, we have ${\cal E}(c\widehat{B}_t^{(\eta)})=\exp\{-\frac{c^2}{2}\eta(t)^2+c\widehat{B}_t^{(\eta)}\}$.
\begin{definition}\label{D:gam}
Let $c> 0$, and suppose $\eta$ is a modulus of continuity, for which the conditions in Theorem \ref{L:MV} hold. Suppose also that $\eta$ satisfies Fernique's condition (see (\ref{E:f30})). Let $\widehat{B}^{(\eta)}$ be the process defined by (\ref{E:ty}). We say that a general time-inhomogeneous Gaussian stochastic volatility model (see Definition \ref{D:gsvm}) belongs to the class ${\cal M}_1^{(\eta)}$ 
if the volatility in the model is described by the process $t\mapsto{\cal E}(c\widehat{B}_t^{(\eta)})$, $t\in[0,T]$.
\end{definition}
\begin{remark}\label{R:tou}
The volatility function in a model from the class ${\cal M}_1^{(\eta)}$ is given by the following formula:
$\sigma(t,u)=\exp\{-\frac{c^2}{2}\eta(t)^2+u\}$, $(t,u)\in[0,T]\times\mathbb{R}$.
This function satisfies the conditions in Assumption C. The previous statement is straightforward, and we leave its proof as an exercise for the interested reader. A model from the class ${\cal M}_1^{(\eta)}$ is super rough if the modulus of continuity $\eta$ satisfies the condition in Lemma \ref{L:mmm}.
\end{remark}
\begin{remark}\label{R:ree}
Assumption A is satisfied for a model from the class ${\cal M}_1^{(\eta)}$. This fact follows from (\ref{E:iso1}), (\ref{E:iso2}),
and part (i) of Theorem \ref{L:MV}.
\end{remark}
\begin{remark}\label{R:ewe}
Using the expression for the covariance function of fractional Brownian motion $B^H$ (see Section \ref{S:aux}), we see that the canonical metric function associated with the process $B^H$ is given by $\delta(t,s)=|t-s|^H$. Define a modulus of continuity on $[0,T]$ by
$\eta(t)=t^H$. Then $\delta(t,s)=\eta(|t-s|)$, and the corresponding process $\widehat{B}^{(\eta)}$ in (\ref{E:ty}) is a constant multiple of the Riemann-Liouville fractional Brownian motion.
\end{remark}
\begin{remark}\label{R:RB}
The models described in Definition \ref{D:gam} have some resemblance to the rough Bergomi model introduced in \cite{BFG}. The latter model is a driftless Gaussian stochastic volatility model, in which the volatility is described by the Wick exponential of the process $cR^H$, where 
$0< H<\frac{1}{2}$, $c> 0$, and $R^H$ is the Riemann-Liouville fractional Brownian motion. The volatility function in the rough Bergomi model is given by $\sigma(t,u)=\exp\{-\frac{c^2t^{2H}}{4H\Gamma(H+\frac{1}{2})^2}+cu\}$, $(t,u)\in[0,T]\times\mathbb{R}$. Here we use the formula for the variance of the process $R^H$ provided in Section \ref{S:aux}.
The previous model is called rough since the sample paths of the volatility process $R^H$ are more rough than those of standard Brownian motion, for which $H=\frac{1}{2}$ 
(see Remark \ref{R:ku}). However, the rough Bergomi model is not super rough. 
\end{remark}

Fix $\beta> 0$ such that $T< e^{-\beta-1}$, and set 
\begin{equation}
\eta_{\beta}(x)=(-\log x)^{-\frac{\beta}{2}},\quad 0< x< T.
\label{E:logg}
\end{equation} 
Then, the modulus of continuity $\eta_{\beta}$ satisfies all the conditions in 
Theorem \ref{L:MV}. The kernel $\tau_{\beta}$ is determined from the equality
$\tau_{\beta}^2(x)=\beta x^{-1}(-\log x)^{-\beta-1}$, $0< x< T$ (see (\ref{E:tr})), while the corresponding Gaussian process in (\ref{E:ty})  is given by
\begin{equation}
\widehat{B}_t^{(\beta)}=\sqrt{\beta}\int_0^t(t-s)^{-\frac{1}{2}}(-\log(t-s))^{-\frac{\beta+1}{2}}dB_s,\,\,0\le t\le T.
\label{E:logh}
\end{equation}
\begin{remark}\label{R:lmb}
The logarithmic modulus of continuity $\eta_{\beta}$ was used in \cite{MV} as an example illustrating Theorem \ref{L:MV}, and the process in
(\ref{E:logh}) was called in \cite{MV} the logarithmic Brownian motion (logBm) with parameter $\beta$. It is not hard to see that if $\beta> 1$, then the conditions in Theorem \ref{L:MV} and Fernique's condition (see (\ref{E:f30})) holds, while if $0<\beta\le 1$, then the conditions in Theorem \ref{L:MV} hold, while Fernique's condition does not hold. The process logBm with $\beta> 1$ is a continuous Volterra Gaussian process (see Lemma \ref{L:ferni}). In addition, Lemma \ref{L:mmm} implies that logBm with $\beta> 1$ is a super rough process. Here we take into account that the function $\eta_{\beta}$ grows near zero faster than any positive power. 
\end{remark}
\begin{remark}\label{R:D}
A time-inhomogeneous Gaussian stochastic volatility model, in which the volatility process $\widehat{B}$ is the process logBm with $\beta> 1$,
will be called a logarithmic model. Such a model is super rough. For example, suppose the volatility in a time-homogeneous Gaussian stochastic volatility model is described by the process 
$t\mapsto\exp\{-\frac{c^2}{2}\left(-\log t\right)^{-\beta}+c\widehat{B}^{(\beta)}_t\},\quad t\in[0,T]$,
where $c> 0$, $\beta> 1$, and the condition $T<e^{-\beta-1}$ holds. Then, we obtain a logarithmic model from the class 
${\cal M}_1^{(\eta_{\beta})}$, where $\eta_{\beta}$ is defined in (\ref{E:logg}).
\end{remark}

We will next introduce one more special class of Gaussian stochastic volatility models.
\begin{definition}\label{D:SaS}
Let $c> 0$, and suppose $\eta$ is a modulus of continuity satisfying the conditions in Theorem \ref{L:MV} and also Fernique's condition 
(see (\ref{E:f30})). We say that a general time-inhomogeneous Gaussian stochastic volatility model (see Definition \ref{D:gsvm}) belongs to the class ${\cal M}_2^{(\eta)}$ if the volatility function is given by $\sigma(u)=c|u|$, $u\in\mathbb{R}$,
and the process $\widehat{B}^{(\eta)}$ defined by (\ref{E:ty}) is the volatility process. 
\end{definition}
\begin{remark}\label{R:lju}
If the modulus of continuity $\eta$ satisfies the conditions in Lemma \ref{L:mmm} and also Fernique's condition, then any model from the 
class ${\cal M}_2^{(\eta)}$ is super rough. For example, let $\eta=\eta_{\beta}$ with $\beta> 1$ (see (\ref{E:logg})). Then, the process 
$\widehat{B}^{(\eta_{\beta})}$ is logBm with parameter $\beta$, and the corresponding logarithmic model from the class 
${\cal M}_2^{(\eta_{\beta})}$ is super rough.
\end{remark}

The volatility function in models described in Definition \ref{D:SaS} is the same as in the Stein-Stein model introduced in 
\cite{SS}, or in the Gaussian models studied in \cite{GVZ1,GVZ2}. The Stein-Stein model is one of the classical stochastic volatility 
models of mathematical finance. It is uncorrelated
($\rho=0$), the drift function in it is constant, and the volatility process is the Ornstein-Uhlenbeck process. Actually, it was claimed by the authors of \cite{SS} that the volatility process in their model is the Ornstein-Uhlenbeck process instantaneously reflected at zero. In fact, the volatility in the Stein-Stein model is the absolute value of the Ornstein-Uhlenbeck process. The previous observation is due to Ball and Roma (see \cite{BR}). The Stein and Stein model with $\rho\neq 0$ was studied in \cite{SZ}. More information about various versions of the Stein-Stein model can be found in \cite{G3}. 
\section{Large Deviation Principles in General Gaussian Stochastic Volatility Models}\label{S:ldr}
The main result obtained in the present section is Theorem \ref{T:27}. It provides a sample path LDP for 
the log-price process $X^{(\varepsilon)}$ in a general time-inhomogeneous Gaussian 
stochastic volatility model (see Definition \ref{D:gsvm}). It follows that Theorem \ref{T:27} can be applied to super rough models, e.g., to logarithmic Gaussian models with $\beta> 1$.  
Several sample path and small-noise LDPs will be derived from 
Theorem \ref{T:27} (see Theorems \ref{T:2}, \ref{T:17}, and \ref{T:1}). In Theorem \ref{T:2}, under an additional assumption that the volatility function is strictly positive, we provide an alternative more simple representation of the rate function in Theorem \ref{T:27}. In Theorems \ref{T:17} and \ref{T:1}, small-noise LDPs for the process $X^{(\varepsilon)}_T$ are obtained. 

Let us denote by $\mathbb{C}_0$ the space of continuous functions on the interval $[0,T]$ equipped with the norm 
$
\displaystyle{||f||_{\mathbb{C}_0}=\max_{t\in[0,T]}|f(t)|}.
$ 
In the sequel, the symbol $\mathbb{H}^1_0$ will stand for the Cameron-Martin space for Brownian motion consisting of absolutely continuous functions $f$ on $[0,T]$ such that $f(0)=0$ and $\dot{f}\in L^2[0,T]$, where $\dot{f}$ is the derivative of the function $f$. 
For a function $f\in\mathbb{H}^1_0$, its norm is defined by 
$
||f||_{\mathbb{H}_0^1}=\{\int_0^T\dot{f}(t)^2dt\}^{\frac{1}{2}}.
$
The following notation will be used throughout the paper:
$$
\widehat{f}(s)=\int_0^sK(s,u)\dot{f}(u)du,\quad s\in[0,T].
$$
\begin{remark}\label{R:remn}
Suppose Assumption A holds. Then it is not hard to see that the linear operator ${\cal K}:L^2[0,T]\mapsto\mathbb{C}_0$ defined by
$({\cal K}g)(t)=\int_0^tK(t,s)g(s)ds$ is compact. Therefore, the image of any closed ball in $\mathbb{H}_0^1$ centered at the zero function and of radius $r> 0$ is precompact in $\mathbb{C}$.
\end{remark}

Consider a measurable functional $\Phi:\mathbb{C}_0^3\mapsto\mathbb{C}_0$ 
defined as follows: For 
$l\in\mathbb{H}^1_0$ and $(f,h)\in\mathbb{C}_0^2$ such that $f\in\mathbb{H}^1_0$ and $h=\widehat{f}$, 
\begin{equation}
\Phi(l,f,h)(t)=\int_0^tb(s,\widehat{f}(s))ds+\bar{\rho}\int_0^t\sigma(s,\widehat{f}(s))\dot{l}(s)ds
+\rho\int_0^t\sigma(s,\widehat{f}(s))\dot{f}(s)ds,
\label{E:bel}
\end{equation}
where $0\le t\le T$. In addition, for all the remaining triples $(l,f,h)$, we set $\Phi(l,f,h)(t)=0$, $t\in[0,T]$. 
Let $g\in\mathbb{C}_0$, and define
\begin{align}
&\widehat{Q}_T(g)=\inf_{l,f\in\mathbb{H}_0^1}\left[\frac{1}{2}\left(\int_0^T\dot{l}(s)^2ds
+\int_0^T\dot{f}(s)^2ds\right):\Phi(l,f,\widehat{f})(t)=g(t),\,
t\in[0,T]\right],
\label{E:vunzi}
\end{align}
if the equation appearing on the right-hand side of (\ref{E:vunzi}) is solvable for $l$ and $f$. If there is no solution, then we set 
$\widehat{Q}_T(g)=\infty$.

The next assertion contains a sample path large deviation principle for the log-price process in a general time-inhomogeneous Gaussian stochastic volatility model. Recall that a rate function on a topological space ${\cal X}$ is a lower semi-continuous mapping 
$I:{\cal X}\mapsto[0,\infty]$. It is assumed that $I$ is not identically infinite. A rate function $I$ is called a good rate function if for every $y\in[0,\infty)$ the level set $L_y=\{x\in{\cal X}:I(x)\le y\}$ is a compact 
subset of ${\cal X}$.
\begin{theorem}\label{T:27}
For any model in Definition \ref{D:gsvm}, the process $\varepsilon\mapsto X^{(\varepsilon)}-x_0$ satisfies the sample path large deviation principle with speed $\varepsilon^{-1}$ 
and good rate function $\widehat{Q}_T$ given by (\ref{E:vunzi}).
The validity of the large deviation principle means that
for every Borel measurable subset ${\cal A}$ of $\mathbb{C}_0$, the following estimates hold:
\begin{align*}
&-\inf_{g\in{\cal A}^{\circ}}\widehat{Q}_T(g)\le\liminf_{\varepsilon\downarrow 0}\varepsilon\log\mathbb{P}\left(
X^{(\varepsilon)}-x_0\in{\cal A}\right) 
 \\
&\le\limsup_{\varepsilon\downarrow 0}\varepsilon\log\mathbb{P}\left(X^{(\varepsilon)}-x_0\in{\cal A}\right)
\le-\inf_{g\in\bar{{\cal A}}}\widehat{Q}_T(g).
\end{align*}
The symbols ${\cal A}^{\circ}$ and $\bar{{\cal A}}$ in the previous estimates stand for the interior and the closure of the set 
${\cal A}$, respectively.
\end{theorem}

The proof of Theorem \ref{T:27} is relegated to Section \ref{S:101}. It is long and rather involved. Next, we will briefly explain the structure of the proof and outline the key ideas on which the proof is based. Our objective is to establish a sample path LDP for the scaled log-price process 
$X^{(\varepsilon)}$ defined at the end of Section \ref{S:aux}. By removing one of the drift terms in the expression for $X^{(\varepsilon)}$ we obtain the process $\widehat{X}^{(\varepsilon)}$ (see (\ref{E:intl})). It is observed in the beginning of Subsection \ref{SS:conti} that the processes 
$X^{(\varepsilon)}$ and $\widehat{X}^{(\varepsilon)}$ satisfy the same sample path LDP. In formula (\ref{E:intl}), the process $\widehat{X}^{(\varepsilon)}$ is represented in terms of the following three processes: independent standard Brownian motions $W$ and $B$, and the volatility process $\widehat{B}$ that depends on $B$. This fact provides sufficient motivation to study the three-component random vector $(W,B,\widehat{B})$. We will prove that this vector is a centered Gaussian random vector in the space $\mathbb{C}^3$ (see Theorem \ref{T:hg}). We will also construct an abstract Wiener space associated with this vector 
(see Theorem \ref{T:aVs}). Next, using a known LDP for abstract Wiener spaces 
formulated in Theorem \ref{T:lpd1}, we establish an LDP for the random vector $(W,B,\widehat{B})$ 
(see Theorem \ref{T:ladp}). It is tempting to try to apply the contraction principle in order to pass from the LDP for 
the random vector $(W,B,\widehat{B})$ to the sample path LDP for the process $\widehat{X}^{(\varepsilon)}$, but it is not clear whether this is possible. Instead, the extended contraction principle (see Theorem 4.2.23 in \cite{DZ}) will be employed. The measurable functional $\Phi:\mathbb{C}_0^3\mapsto\mathbb{C}_0$ given by (\ref{E:bel}) imitates the structure of the log-price process $X$. We will use a sequence  
$\Phi_m:\mathbb{C}_0^3\mapsto\mathbb{C}_0$, $m\ge 2$, of discrete approximations of the functional $\Phi$ (see (\ref{E:lopi}) for the definition of this sequence). It is explained in Lemma \ref{L:ecp} in what sense the sequence $\{\Phi_m\}$ of continuous functionals approximates the functional $\Phi$. This lemma shows that one of the conditions in the contraction principle holds in our environment. A crucial part of the proof of Theorem \ref{T:27} is to establish that the sequence of the processes 
$\varepsilon\mapsto\Phi_m\left(\sqrt{\varepsilon}W,\sqrt{\varepsilon}B,
\sqrt{\varepsilon}\widehat{B}\right)$ with state space $\mathbb{C}_0$ is an exponentially good approximation to the process $\varepsilon\mapsto\widehat{X}^{(\varepsilon)}$. This is done in Lemma \ref{L:borr}. Finally, we can use the extended contraction principle to complete the proof of Theorem \ref{T:27} (see the end of Subsection \ref{SS:conti}).
\begin{remark}\label{R:yui}
The fact that the LDP in Theorem \ref{T:27} applies to super rough models can be regarded as the main novelty in the present paper. Similar LDPs obtained in the earlier papers \cite{FZ,G1,G2,CP} cannot be applied to such models (see the discussion in the introduction).  Although there are certain resemblances between the proof of Theorem \ref{T:27} and the proofs of Theorem 4.5 in 
\cite{FZ} and Theorem 2.9 in \cite{G2}, there are also important differences. The key result illustrating these differences is 
Corollary \ref{C:corol} (see Subsection \ref{SS:ert}). 
In this corollary, a large deviation style formula is established for a certain maximal function of the increments of the 
volatility process $\widehat{B}$. Corollary \ref{C:corol} is general enough to be applied to all the super rough volatility processes of our interest in the present paper. This corollary plays an important part in the proof of Theorem \ref{T:27}. For instance, 
Corollary \ref{C:corol} is used in Subsection \ref{SS:conti} to establish formulas (\ref{E:tru}) -- (\ref{E:truss}). These formulas are significant components of the proof of Lemma \ref{L:borr} in Subsection \ref{SS:conti}. We derive Corollary \ref{C:corol} 
from the sample path LDP for the process $\varepsilon\mapsto\sqrt{\varepsilon}\widehat{B}$ that is formulated in Remark \ref{R:foll}. This is a new approach to proving maximal estimates for the increments of the volatility process. Note that such estimates were also used in 
\cite {FZ,G1,G2} under stronger restrictions. In \cite{FZ}, the volatility is modeled by fractional Brownian motion, which is a process with stationary increments. The previous fact makes the proof of maximal estimates considerably simpler (see (B-6) in \cite{FZ}). The model studied in \cite{FZ} is not super rough. This has already been mentioned in the introduction. In our earlier papers \cite{G1} and \cite{G2}, a maximal estimate for the increments of the volatility process was also employed (see Lemma 23 in \cite{G1}). In the proof of the previous lemma, we used the maximal estimates for the increments of certain Banach space valued stochastic processes obtained in \cite{CC}, Lemma 2.2. However, Lemma 23 in \cite{G1} is not general enough to be applied to super rough models appearing in the present paper. 
\end{remark}
\begin{remark}\label{R:inh}
Unlike the earlier models studied in \cite{FZ,G1,G2,CP} the models considered in the present paper are time-inhomogeneous. This does not cause serious problems in the proof of the LDP in Theorem \ref{T:27}. However, various technical difficulties arise, especially in the estimates used in the proof of Lemma \ref{L:borr} in Subsection \ref{SS:conti}.
\end{remark}
The next assertion can be derived from Theorem \ref{T:27}. This derivation will be explained in the sequel.
\begin{theorem}\label{T:2}
Suppose that the volatility function $\sigma$ in a model satisfying the conditions in Definition \ref{D:gsvm} is strictly positive on $[0,T]\times\mathbb{R}$. 
Then the process $\varepsilon\mapsto X^{(\varepsilon)}-x_0$ satisfies the sample path large deviation principle 
with speed $\varepsilon^{-1}$ 
and good rate function $Q_T$ given by $Q_T(g)=\infty$, for all $g\in\mathbb{C}_0\backslash\mathbb{H}_0^1$, and 
\begin{align}
&Q_T(g)=\inf_{f\in\mathbb{H}_0^1}\left[\frac{1}{2}\int_0^T\left[\frac{\dot{g}(s)-b(s,\widehat{f}(s))-\rho\sigma(s,\widehat{f}(s))\dot{f}(s)}
{\bar{\rho}\sigma(s,\widehat{f}(s))}\right]^2ds
+\frac{1}{2}\int_0^T\dot{f}(s)^2ds\right],
\label{E:vunzil}
\end{align}
for all $g\in\mathbb{H}_0^1$. 
\end{theorem}

Our next goal is to formulate small-noise large deviation principles for the log-price process in models introduced in Definition 
\ref{D:gsvm} (see Theorems \ref{T:17} and \ref{T:1} below). Suppose the volatility function is not identically zero on 
$[0,T]\times\mathbb{R}$. Note that Assumption C includes the validity of the previous condition. Recall that throughout the paper we assume that $K$ is an admissible kernel.
\begin{definition}\label{D:loi}
Denote by $L_1$ the set of all functions $f\in\mathbb{H}_0^1$ satisfying the following condition:
$\sigma(s,\widehat{f}(s))=0$, $s\in[0,T]$,
and set $L_2=\mathbb{H}_0^1\backslash L_1$.
\end{definition}
By the continuity of the functions $\sigma$ and $\widehat{f}$, the condition in Definition \ref{D:loi} is equivalent to 
the following equality:
$
\int_0^T\sigma(s,\widehat{f}(s))^2ds=0.
$

The next two lemmas concern the sets $L_1$ and $L_2$ introduced in Definition \ref{D:loi}. The proofs of these lemmas are straightforward, and we omit them.
\begin{lemma}\label{L:short}
The set $L_2$ is not empty. 
\end{lemma}

Unlike the set $L_2$, the set $L_1$ may be empty. For example, this happens if the function $\sigma$ 
is strictly positive. On the other hand, if the volatility function $\sigma$ is such that 
\begin{equation}
\sigma(s,0)=0,\quad\mbox{for all}\quad s\in[0,T],
\label{E:sigg}
\end{equation} 
then $f_0\in L_1$, where the symbol $f_0$ stands for the function equal to zero identically. 
\begin{lemma}\label{L:emp}
Suppose the volatility function $\sigma$ is such that condition (\ref{E:sigg}) holds.
Then, there exist an admissible kernel $K$ and a function $f\neq f_0$ such that $f\in L_1$.
\end{lemma}

For $y\in\mathbb{R}$ and $f\in\mathbb{H}_0^1$, set
$$
\Psi(y,f,\widehat{f})=\int_0^T[b(s,\widehat{f}(s))+\rho\sigma(s,\widehat{f}(s))\dot{f}(s)]ds+\bar{\rho}
\left\{\int_0^T\sigma(s,\widehat{f}(s))^2ds\right\}^{\frac{1}{2}}y.
$$
\begin{lemma}\label{L:denote}
For every $x\in\mathbb{R}$ and $f\in L_2$, there exists a number $y\in\mathbb{R}$ such that
$\Psi(y,f,\widehat{f})=x$.
\end{lemma}

\it Proof. \rm For $f\in L_2$, we can solve the equation in Lemma \ref{L:denote} for $y$, which gives
$$
y=\frac{x-\int_0^T[b(s,\widehat{f}(s))+\rho\sigma(s,\widehat{f}(s))\dot{f}(s)]ds}
{\bar{\rho}
\{\int_0^T\sigma(s,\widehat{f}(s))^2ds\}^{\frac{1}{2}}}.
$$
This establishes Lemma \ref{L:denote}.

It follows from Lemma \ref{L:denote} that
\begin{equation}
y^2=\Lambda(x,f),
\label{E:repre}
\end{equation}
where
\begin{equation}
\Lambda(x,f)=\frac{\left(x-\int_0^T[b(s,\widehat{f}(s))+\rho\sigma(s,\widehat{f}(s))\dot{f}(s)]ds\right)^2}
{\bar{\rho}^2
\int_0^T\sigma(s,\widehat{f}(s))^2ds}.
\label{E:retro}
\end{equation}

For every $x\in\mathbb{R}$, set
\begin{align}
&\widehat{I}_T(x)=\inf_{y\in\mathbb{R},f\in\mathbb{H}_0^1}\left[\frac{1}{2}\left(y^2
+\int_0^T\dot{f}(s)^2ds\right):\Psi(y,f,\widehat{f})=x\right].
\label{E:vunzics}
\end{align}
Using Lemma \ref{L:denote}, we see that the function $\widehat{I}_T$ is well-defined.

We will next find an alternative representation for the function $\widehat{I}_T$. 
It is easy to see that the equation 
\begin{equation}
\Psi(y,f,\widehat{f})=x
\label{E:labe}
\end{equation} 
with the restriction $f\in L_1$ can be rewritten as follows:
\begin{equation}
\int_0^Tb(s,\widehat{f}(s))ds=x.
\label{E:trt}
\end{equation}
Let $x\in\mathbb{R}$, and suppose the equation in (\ref{E:labe}) has a solution $f\in L_1$. 
Then, we denote the set of all solutions from $L_1$
by $M_x$.
\begin{lemma}\label{L:dgh}
The following formulas hold for the function $\widehat{I}_T$: 
\begin{align}
&\widehat{I}_T(x)=\frac{1}{2}\min\left\{\inf_{f\in M_x}\int_0^T\dot{f}(s)^2ds,   
\inf_{f\in L_2}\left[\Lambda(x,f)
+\int_0^T\dot{f}(s)^2ds\right]\right\},
\label{E:pous}
\end{align}
if $M_x\neq\emptyset$, and
\begin{align}
&\widehat{I}_T(x)=\frac{1}{2}
\inf_{f\in L_2}\left[\Lambda(x,f)
+\int_0^T\dot{f}(s)^2ds\right],
\label{E:pouu}
\end{align}
if the equation in (\ref{E:labe}) does not have any solutions in $L_1$. In (\ref{E:pous}) and (\ref{E:pouu}), the function $\Lambda$
is defined by (\ref{E:retro}).
\end{lemma}

\it Proof. \rm Lemma \ref{L:dgh} can be established using formulas (\ref{E:repre}) -- (\ref{E:trt}).

Finally, we are ready to formulate small-noise large deviation principles for general time-inhomogeneous Gaussian models.
\begin{theorem}\label{T:17}
For any model in Definition \ref{D:gsvm}, the process $\varepsilon\mapsto X^{(\varepsilon)}_T-x_0$ satisfies the small-noise large deviation principle with speed $\varepsilon^{-1}$ 
and good rate function $\widehat{I}_T$ given by (\ref{E:vunzics}). The validity of such a large deviation principle means that
for every Borel measurable subset $A$ of $\mathbb{R}$, the following estimates hold:
\begin{align*}
&-\inf_{x\in A^{\circ}}\widehat{I}_T(x)\le\liminf_{\varepsilon\downarrow 0}\varepsilon\log\mathbb{P}\left(
X_T^{(\varepsilon)}-x_0\in A\right) 
 \\
&\le\limsup_{\varepsilon\downarrow 0}\varepsilon\log\mathbb{P}\left(X_T^{(\varepsilon)}-x_0\in A\right)
\le-\inf_{x\in\bar{A}}\widehat{I}_T(x).
\end{align*}
The symbols $A^{\circ}$ and $\bar{A}$ in the previous estimates stand for the interior and the closure of the set $A$, respectively.
\end{theorem}
\begin{theorem}\label{T:1}
Suppose that the volatility function $\sigma$ in a model satisfying the conditions in Definition \ref{D:gsvm} is strictly positive on $[0,T]\times\mathbb{R}$. Set
\begin{equation}
I_T(x)=\inf_{f\in\mathbb{H}_0^1}\left[\frac{\left(x-\int_0^T[b(s,\widehat{f}(s))+\rho\sigma(s,\widehat{f}(s))\dot{f}(s)]ds\right)^2}
{2\bar{\rho}^2\int_0^T\sigma(s,\widehat{f}(s))^2ds}+\frac{1}{2}\int_0^T\dot{f}(s)^2ds\right].
\label{E:vunz}
\end{equation}
Then $I_T$ is a good rate function. Moreover, a small-noise large deviation principle with speed $\varepsilon^{-1}$ 
and rate function $I_T$ given by (\ref{E:vunz}) holds for the process
$\varepsilon\mapsto X_T^{(\varepsilon)}-x_0$. 
\end{theorem}
\begin{remark}\label{R:remis}
A set ${\cal A}\subset\mathbb{C}_0$ is called a set of continuity for the rate function $\widehat{Q}_T$ if
\begin{equation}
\inf_{g\in{\cal A}^{\circ}}\widehat{Q}_T(g)=\inf_{g\in\bar{{\cal A}}}\widehat{Q}_T(g).
\label{E:0}
\end{equation}
For such a set, Theorem \ref{T:27} implies that
\begin{equation}
\lim_{\varepsilon\downarrow 0}\varepsilon
\log\mathbb{P}\left(X^{(\varepsilon)}-x_0\in{\cal A}\right)=-\inf_{g\in{\cal A}}\widehat{Q}_T(g).
\label{E:00}
\end{equation}
Similar statements can be derived from Theorems \ref{T:2}, \ref{T:17}, and \ref{T:1} for the sets of continuity for the corresponding rate functions.
\end{remark}

 Our next goal is to show that Theorems \ref{T:2}, \ref{T:17}, and \ref{T:1} follow from Theorem \ref{T:27}.

\it Derivation of Theorem \ref{T:2} from Theorem \ref{T:27}. \rm
Suppose Theorem \ref{T:27} holds, and the volatility
function $\sigma$ is strictly positive on $[0,T]\times\mathbb{R}$. Then, Theorem \ref{T:2} follows from the following lemma.
\begin{lemma}\label{L:lop}
Under the conditions in Theorem \ref{T:2}, $\widehat{Q}_T=Q_T$, where the functions $\widehat{Q}_T$ and $Q_T$ are defined by (\ref{E:vunzi})
and (\ref{E:vunzil}), respectively.
\end{lemma}

\it Proof. \rm
Suppose $g\in\mathbb{C}_0$, $l,f\in\mathbb{H}_0^1$, and $\Phi(l,f,\widehat{f})(t)=g(t)$, for all $t\in[0,T]$. Then, it is easy to see that 
$g\in\mathbb{H}_0^1$. Moreover, if for $g\in\mathbb{H}_0^1$ the previous equality holds, then 
$$
\dot{l}(t)=\frac{\dot{g}(t)-b(t,\widehat{f}(t))-\rho\sigma(t,\widehat{f}(t)))\dot{f}(t)}
{\bar{\rho}\sigma(t,\widehat{f}(t))}
$$
almost everywhere with respect to the Lebesgue measure on $[0,T]$. Here we use the fact that the volatility function is strictly positive. Now, it is clear that the functions $Q_T$ and $\widehat{Q}_T$ are equal.

This completes the proof of Lemma \ref{L:lop}.

\it Derivation of Theorem \ref{T:1} from Theorem \ref{T:17}. \rm Theorem \ref{T:1} follows from Theorem \ref{T:17} and Lemma \ref{L:dgh}. Indeed, if the volatility function is strictly positive, then
$L_2=\mathbb{H}_0^1$. In addition, the equation in (\ref{E:labe}) does not have any solutions in $L_1$.

\it Derivation of Theorem \ref{T:17} from Theorem \ref{T:27}. \rm 
Let us consider a mapping $V:\mathbb{C}_0\mapsto\mathbb{R}$ defined by $V(\varphi)=\varphi(T)$, $\varphi\in\mathbb{C}_0$. 
It is clear that the mapping $V$ is continuous. Suppose Theorem \ref{T:27} holds. Then, applying the contraction principle we see that the process 
$\varepsilon\mapsto X_T^{(\varepsilon)}-x_0$ satisfies the LDP with speed $\varepsilon^{-1}$ and good rate function given by
\begin{align*}
&\widetilde{I}_T(x)=\inf_{\{g\in\mathbb{H}_0^1,g(T)=x\}}\widehat{Q}_T(g) \\
&=\inf_{\{g,l,f\in\mathbb{H}_0^1,g(T)=x\}}\left[\frac{1}{2}\left(\int_0^T\dot{l}(s)^2ds
+\int_0^T\dot{f}(s)^2ds\right)
: \Phi(l,f,\widehat{f})(t)=g(t),\,t\in[0,T]\right],
\end{align*}
where $x\in\mathbb{R}$.

Recall that the function $\widehat{I}_T$ is defined in (\ref{E:vunzics}). We will next prove that for every $x\in\mathbb{R}$, the following equality holds: 
\begin{equation}
\widehat{I}_T(x)=\widetilde{I}_T(x)\quad x\in\mathbb{R}.
\label{E:ryu}
\end{equation} 
It is clear that Theorem \ref{T:17} follows from the previous equality. Fix $x$ and $g$ with $g(T)=x$, and suppose the equation 
\begin{equation}
\Phi(l,f,\widehat{f})(t)=g(t),\,t\in[0,T],
\label{E:v1}
\end{equation}
is not solvable for $(l,f)$. Then, we have $\widetilde{I}_T(x)=\infty$, and hence $\widehat{I}_T(x)\le\widetilde{I}_T(x)$.
Note that for $g\in\mathbb{C}_0\backslash\mathbb{H}_0^1$, the equation in (\ref{E:v1}) is not solvable.

In the sequel, we will use the sets $L_1$ and $L_2$ defined in the proof of Lemma \ref{L:denote}.
Suppose there exists a solution $(l,f)$ to the equation in (\ref{E:v1}) such that $f\in L_1$. Then, we have
$g(t)=\int_0^tb(s,\widehat{f}(s))ds$, for all $t\in[0,T]$.
Moreover, $x=g(T)=\int_0^Tb(s,\widehat{f}(s))ds$. Fix the function $f$ and consider the equation $x=\Psi(y,f,\widehat{f})$. It follows
from the reasoning above that any $y\in\mathbb{R}$ is a solution to the previous equation.

Next, suppose there is a solution $(l,f)$ to the equation in (\ref{E:v1}) with $f\in L_2$. Then we have
\begin{equation}
0<\int_0^T\sigma(s,\widehat{f}(s))^2ds<\infty,
\label{E:d1}
\end{equation}
and for all $t\in[0,T]$, the following equality holds:
\begin{equation}
g(t)=\int_0^tb(s,\widehat{f}(s))ds
+\rho\int_0^t\sigma(s,\widehat{f}(s))\dot{f}(s)ds+\bar{\rho}\int_0^t\sigma(s,\widehat{f}(s))\dot{l}(s)ds.
\label{E:d3}
\end{equation}
Hence,
\begin{equation}
x=g(T)=\int_0^Tb(s,\widehat{f}(s))ds
+\rho\int_0^T\sigma(s,\widehat{f}(s))\dot{f}(s)ds+\bar{\rho}\int_0^T\sigma(s,\widehat{f}(s))\dot{l}(s)ds.
\label{E:d2}
\end{equation}

Consider the equation 
\begin{equation}
x=\Psi(y,f,\widehat{f}),
\label{E:ui1}
\end{equation} 
where $x$ is the same as in (\ref{E:d2}), and $f\in L_2$ is the same as in 
(\ref{E:d3}). We have already established that the equation in (\ref{E:ui1}) with a fixed function $f\in L_2$ is solvable for $y$, 
and we have $y^2=\Lambda(x,f)$ (see (\ref{E:repre})). 

It will be shown next that 
\begin{equation}
y^2\le\int_0^T\dot{l}(s)^2ds.
\label{E:d7}
\end{equation} 
Indeed, the equation in (\ref{E:ui1}) is as follows:
\begin{align}
&x=\int_0^Tb(s,\widehat{f}(s))ds
+\rho\int_0^T\sigma(s,\widehat{f}(s))\dot{f}(s)ds 
+\bar{\rho}
\left\{\int_0^T\sigma(s,\widehat{f}(s))^2ds\right\}^{\frac{1}{2}}y.
\label{E:d4}
\end{align}
Next, using (\ref{E:d2}) and (\ref{E:d4}), we obtain
$$
\left\{\int_0^T\sigma(s,\widehat{f}(s))^2ds\right\}^{\frac{1}{2}}y
=\int_0^T\sigma(s,\widehat{f}(s))\dot{l}(s)ds.
$$
Therefore, 
$$
y^2\int_0^T\sigma(s,\widehat{f}(s))^2ds
\le\left(\int_0^T\sigma(s,\widehat{f}(s))^2ds\right)\left(\int_0^T\dot{l}(s)^2ds\right).
$$
Next, using (\ref{E:d1}), we obtain (\ref{E:d7}).

Finally, by taking into account the reasoning above, we obtain the following estimate:
\begin{equation}
\widehat{I}_T(x)\le\widetilde{I}_T(x),\quad x\in\mathbb{R}.
\label{E:opp1}
\end{equation}

We will next prove the opposite inequality. Let $x\in G$, and consider the equation in (\ref{E:d4}). If there is no solution $(y,f)$
to this equation, then $\widehat{I}_T(x)=\infty$, and it follows that $\widetilde{I}_T(x)\le\widehat{I}_T(x)$.

Next, suppose there exists a solution $(y,f)$ to the equation in (\ref{E:d4}) with $f\in L_1$. Then we have
$x=\int_0^Tb(s,\widehat{f}(s))ds$.
Set $g(t)=\int_0^tb(s,\widehat{f}(s))ds$, $t\in[0,T]$. Then $g\in\mathbb{H}^1_0$, and $g(T)=x$. It is not hard to see that
(\ref{E:d3}) holds for the function $g$ with the same function $f$ and any function $l\in\mathbb{H}^1_0$.

Next, suppose $(y,f)$ is a solution to the equation in (\ref{E:d4}) with $f\in L_2$. Then, (\ref{E:d1})
holds, and there exists a function $l\in\mathbb{H}^1_0$ such that for all $s\in[0,T]$,
\begin{equation}
\dot{l}(s)=\frac{y}
{\left\{\int_0^T\sigma(u,\widehat{f}(u))^2du\right\}^{\frac{1}{2}}}\sigma(s,\widehat{f}(s)).
\label{E:rer}
\end{equation}
Fix such a function $l$, and set
\begin{equation}
g(t)=\int_0^tb(s,\widehat{f}(s))ds
+\rho\int_0^t\sigma(s,\widehat{f}(s))\dot{f}(s)ds+\bar{\rho}\int_0^t\sigma(s,\widehat{f}(s))\dot{l}(s)ds,
\label{E:bri}
\end{equation}
for all $t\in[0,T]$. It is clear that the pair $(l,f)$ solves the equation in (\ref{E:v1}) for the function $g$.
Moreover, $\int_0^T\dot{l}(s)^2ds=y^2$. It remains to prove that $g(T)=x$. 
By (\ref{E:bri}), the previous equality is equivalent to
\begin{equation}
x=\int_0^Tb(s,\widehat{f}(s))ds
+\rho\int_0^T\sigma(s,\widehat{f}(s))\dot{f}(s)ds+\bar{\rho}\int_0^T\sigma(s,\widehat{f}(s))\dot{l}(s)ds.
\label{E:labb}
\end{equation}
Next, using (\ref{E:rer}), we can rewrite the equality in (\ref{E:labb}) as follows:
\begin{align}
x&=\int_0^Tb(s,\widehat{f}(s))ds
+\rho\int_0^T\sigma(s,\widehat{f}(s))\dot{f}(s)ds+\bar{\rho}
\left\{\int_0^T\sigma(s,\widehat{f}(s))^2ds\right\}^{\frac{1}{2}}y,
\label{E:brius}
\end{align}
for all $t\in[0,T]$. The equality in (\ref{E:brius}) holds, since the pair $(y,f)$ solves the equation in (\ref{E:d4}).
It follows that the equality $g(T)=x$ holds true.

Summarizing what was said above, we see that
$\widetilde{I}_T(x)\le\widehat{I}_T(x)$, $x\in\mathbb{R}$.
Next, combining (\ref{E:opp1}) and the previous inequality, we obtain (\ref{E:ryu}). As we have already mentioned above, the equality
in (\ref{E:ryu}) implies Theorem \ref{T:17}.
\begin{remark}\label{R:rio}
Since the volatility functions in the formulations of Theorems \ref{T:2} and \ref{T:1} are strictly positive, typical examples here are exponential functions. On the other hand, Theorems \ref{T:27} 
and \ref{T:17} allow for
power type volatility functions, e.g., the function $\sigma(t,u)=|u|$, $(t,u)\in[0,T]\times\mathbb{R}$.
\end{remark}

\begin{remark}\label{R:lower}
The rate functions $\widehat{I}_T$ and $I_T$ in Theorems \ref{T:17} and \ref{T:1} are lower semicontinuous on $\mathbb{R}$. 
The rate function $I_T$ is also upper semicontinuous since it is the greatest lower bound of a family of continuous functions. 
It follows that $I_T$ is a continuous function. On the other hand, it is not clear whether the function
$\widehat{I}_T$ is always continuous. In the next lemma, we will show that in the case, where the drift coefficient is constant, the function
$\widehat{I}_T$ may be discontinuous only at one point.
\end{remark}

\begin{lemma}\label{L:ptu}
Suppose $b(s,u)=r$, for all $(s,u)\in[0,T]\times\mathbb{R}$ and some $r\in\mathbb{R}$. Then, the good rate function $\widehat{I}_T$ 
is continuous on $\mathbb{R}\backslash\{rT\}$.
\end{lemma}

\it Proof. \rm Since the function $\widehat{I}_T$ is lower semicontinuous on $\mathbb{R}$, we only need to establish that it is upper semicontinuous on $\mathbb{R}\backslash\{rT\}$.
It is not hard to show that, under the conditions in Lemma \ref{L:ptu}, a solution $f\in L_1$ to the equation in (\ref{E:labe}) exists only when $x=rT$. Using Lemma \ref{L:dgh}, we see that for $x\neq rT$,
$\widehat{I}_T(x)=\frac{1}{2}\inf_{f\in L_2}\left[\Lambda(x,f)
+\int_0^T\dot{f}(s)^2ds\right]$.
The function on the right-hand side of the previous equality is nonnegative and upper semicontinuous on $\mathbb{R}$ being the infimum of a family of continuous functions. Therefore, the function $\widehat{I}_T$ is upper semicontinuous on $\mathbb{R}\backslash\{rT\}$.

The proof of Lemma \ref{L:ptu} is thus completed.

\section{Applications}\label{S:addik}
Our first goal in the present section is to study a small-noise asymptotic behavior of the initial price of a scaled barrier option. We consider four standard barrier options (up-and-in, up-and-out, down-and-in, down-and-out). For every fixed $\varepsilon\in(0,1]$, the process
$S^{(\varepsilon)}$ defined at the end of Section \ref{S:aux} will play the role of the underlying asset price process.
It will be assumed in this section that the drift coefficient $b$ in the model  
satisfies $b(s,u)=r$, for all $(s,u)\in[0,T]\times\mathbb{R}$, where $r\ge 0$ is the interest rate. Let us set the barrier at $K> 0$, 
and let $T> 0$ be the maturity of the option. Suppose the following inequality holds: $s_0< K$. The up-and-in binary barrier option pays a fixed amount $G$ of cash if the underlying process touches the barrier at some time during the life of the option. For every $\varepsilon\in(0,T]$, the price of such an option at the initial moment $t=0$ is given by the formula 
$V_1(\varepsilon)=Ge^{-rT}\mathbb{P}(\max_{t\in[0,T]}S_t^{(\varepsilon)}\ge K)$. Indeed, the above-mentioned price is equal to the expectation of the discounted pay-off of the option, and for the up-and-in binary barrier option, this pay-off equals 
$G\mathbb{1}_{\{\max_{t\in[0,T]}S_t^{(\varepsilon)}\ge K\}}$. The up-and-out binary barrier option pays a fixed amount $G$ of cash if the asset price process never touches the barrier during the life of the option. The price of the option in this case is given by 
$V_2(\varepsilon)=Ge^{-rT}\mathbb{P}(\max_{t\in[0,T]}S_t^{(\varepsilon)}<K)$.  
Now, let $K< s_0$. Then, the down-and-in and down-and-out binary options are defined similarly to the definitions above. More precisely,
the prices of the down-and-in and down-and-out options are given by
$V_3(\varepsilon)=Ge^{-rT}\mathbb{P}(\min_{t\in[0,T]}S_t^{(\varepsilon)}\le K)$ and
$V_4(\varepsilon)=Ge^{-rT}\mathbb{P}(\min_{t\in[0,T]}S_t^{(\varepsilon)}> K)$, respectively. 
Let $s_0< K$, and denote: 
$$
{\cal A}_T^{(1)}=\left\{f\in\mathbb{C}_0:f(s)+x_0=\log K\,\,\mbox{for some}\,\,s\in(0,T]\right\}
$$
and
$$
{\cal A}_T^{(2)}=\left\{f\in\mathbb{C}_0:f(s)+x_0<\log K
\,\,\mbox{for all}\,\,s\in(0,T]\right\}.
$$
Similarly, in the case where $s_0> K$, the following sets will be used:
$$
{\cal A}_T^{(3)}=\left\{f\in\mathbb{C}_0:f(s)+x_0=\log K
\,\,\mbox{for some}\,\,s\in(0,T]\right\}
$$
and 
$$
{\cal A}_T^{(4)}=\left\{f\in\mathbb{C}_0:f(s)+x_0>\log K
\,\,\mbox{for all}\,\,s\in(0,T]\right\}.
$$
It is not hard to see that the sets ${\cal A}_1$ and ${\cal A}_3$ are closed in the space $\mathbb{C}_0$, while
the sets ${\cal A}_2$ and ${\cal A}_4$ are open.

For every $1\le k\le 4$, set $L_k=\inf_{f\in{\cal A}_T^{(k)}}Q_T(f)$. It is clear that $L_k\ge 0$.
The next assertion provides large deviation style formulas for binary barrier options. 
\begin{theorem}\label{T:333}
Under the conditions in Theorem \ref{T:2} and the restrictions in the definitions of binary digital options ($s_0< K$, or $K< s_0$), 
the following formula holds:
\begin{equation}
\lim_{\varepsilon\rightarrow 0}\varepsilon\log V_k(\varepsilon)=-L_k,\quad 1\le j\le 4,
\label{E:chetyre}
\end{equation}
where $Q_T$ is the rate function given by (\ref{E:vunzil}).
\end{theorem}

\it Proof of Theorem \ref{T:333}. \rm We will prove formula (\ref{E:chetyre}) for the up-and-in option, that is, in the case where $k=1$. 
In the remaining cases, the proof is similar. The next lemma will be used in the proof.
\begin{lemma}\label{L:concon}
The functional $Q_T:\mathbb{H}^1_0\mapsto\mathbb{R}$ is continuous.
\end{lemma}

\it Proof of Lemma \ref{L:concon}. \rm The lower semi-continuity of the functional $Q_T$ in Lemma \ref{L:concon} follows from the fact 
that $Q_T$ 
is a rate function on $\mathbb{C}_0$ and the continuous embedding $\mathbb{H}^1_0\subset\mathbb{C}_0$. 

We will next prove the upper semi-continuity of $Q_T$ on $\mathbb{H}^1_0$. For every $f\in\mathbb{H}^1_0$, define the functional
${\cal D}_f:\mathbb{H}^1_0\mapsto\mathbb{R}$ by
$$
{\cal D}_f(g)=\int_0^T\left[\frac{\dot{g}(s)-b(s,\widehat{f}(s))-\rho\sigma(s,\widehat{f}(s))\dot{f}(s)}
{\bar{\rho}\sigma(s,\widehat{f}(s))}\right]^2ds
+\int_0^T\dot{f}(s)^2ds.
$$
It is not hard to see that in order to complete the proof of Lemma \ref{L:concon}, it suffices to establish that for every 
$f\in\mathbb{H}^1_0$,
the functional ${\cal D}_f$ is continuous. The function $\widehat{f}$ is a continuous function on $[0,T]$. Therefore
there exist $\delta> 0$ and $M> 0$ such that
$b(s,\widehat{f}(s))< M$ and $\delta<\sigma(s,\widehat{f}(s))< M$, $s\in[0,T]$.

Suppose $g_k\rightarrow g$ in $\mathbb{H}^1_0$. Then we have
\begin{align*}
&|{\cal D}_f(g)-{\cal D}_f(g_k)| \\
&\le\frac{1}{\delta^2(1-\rho^2)}\int_0^T|\dot{g}(s)-\dot{g}_k(s)|
|\dot{g}(s)+\dot{g}_k(s)-2b(s,\widehat{f}(s))-2\rho\sigma(s,\widehat{f}(s))\dot{f}(s)|ds \\
&\le\frac{1}{\delta^2(1-\rho^2)}||g-g_k||_{\mathbb{H}^1_0} \\
&\times\left(||g||_{\mathbb{H}^1_0}+\sup_k||g_k||_{\mathbb{H}^1_0}
+2\left\{\int_0^Tb(s,\widehat{f}(s))^2ds\right\}^{\frac{1}{2}}
+2\left\{\int_0^T\sigma(s,\widehat{f}(s))^2\dot{f}(s)^2ds\right\}^{\frac{1}{2}}\right) \\
&\le\frac{1}{\delta^2(1-\rho^2)}||g-g_k||_{\mathbb{H}^1_0}\left(||g||_{\mathbb{H}^1_0}+\sup_k||g_k||_{\mathbb{H}^1_0}
+2M+2M||f||_{\mathbb{H}^1_0}\right).
\end{align*}
Now, it is clear that ${\cal D}_f(g_k)\rightarrow{\cal D}_f(g)$ as $k\rightarrow\infty$, and hence the functional ${\cal D}_f$
is continuous on the space $\mathbb{H}^1_0$. It follows that the functional $Q_T$ is upper semi-continuous since it can 
be represented as the infimum of a family of continuous on $\mathbb{H}^1_0$ functionals.

This completes the proof of Lemma \ref{L:concon}.

Let us return to the proof of Theorem \ref{T:333}. It is easy to see that 
\begin{equation}
V_1(\varepsilon)=Ge^{-rT}\mathbb{P}(X^{(\varepsilon)}_{\cdot}-x_0\in{\cal A}_T^{(1)}).
\label{E:lip}
\end{equation}
We will next show that ${\cal A}_T^{(1)}$ is a set of continuity for the rate function $Q_T$ (see the definition in Remark \ref{R:remis}). 
It has already been mentioned that ${\cal A}_T^{(1)}$ is a closed subset of the space $\mathbb{C}_0$. Moreover, the interior 
$[{\cal A}^{(1)}_T]^{\circ}$ of the set ${\cal A}_T$ 
consists of all $f\in\mathbb{C}_0$, for which there exists $s< T$ such that $f(s)+x_0>\log K$. In addition, 
a function $f\in\mathbb{C}_0$ belongs to the boundary of the set ${\cal A}_T^{(1)}$ 
if and only if the function $f+x_0$ hits the point $\log K$
before $T$ or at $s=T$, 
but never exits the set $(-\infty,\log K]$. It is not hard to see that the set 
$[{\cal A}^{(1)}_T]^{\circ}\cap\mathbb{H}_0^1$ is dense in the set
${\cal A}_T^{(1)}\cap\mathbb{H}_0^1$ in the topology of the space $\mathbb{H}_0^1$. Now, using Lemma \ref{L:concon} and taking into account the discussion in Remark
\ref{R:remis}, it is easy to prove 
that the equality in (\ref{E:0}) holds for the set ${\cal A}_T^{(1)}$ and the rate function $Q_T$. It follows that ${\cal A}_T^{(1)}$ is a set of continuity for the rate function
$Q_T$. Therefore, the equality in (\ref{E:00}) 
is valid for the set ${\cal A}_T^{(1)}$. Finally, the equality in (\ref{E:chetyre}) for $k=1$ 
can be established using (\ref{E:lip}) and Theorem \ref{T:2}.

The price of a binary up-and-in barrier option is related to the exit time probability function of the  
log-price process $X^{(\varepsilon)}$ from the set $U=(-\infty,\log K)$. Since $s_0\in(0,K)$, we have $x_0\in U$.
The exit time of the process $X^{(\varepsilon)}$ from the set $U$ is defined by
$$
\tau^{(\varepsilon)}=\inf\left\{s\in(0,T]:X_s^{(\varepsilon)}\notin U\right\}\quad\varepsilon\in(0,1].
$$
Here we assume that $\inf\emptyset=\infty$. The exit time probability function associated with the process $X^{(\varepsilon)}$ is given by 
$v_{\varepsilon}(t)=\mathbb{P}(\tau^{\varepsilon}\le t)$, $t\in(0,T]$.
It is not hard to see that
\begin{equation}
V_1(\varepsilon)=Ge^{-rT}v_{\varepsilon}(T)\quad\varepsilon\in(0,1].
\label{E:exit} 
\end{equation}

The next assertion provides a large deviation style formula in the small-noise regime for the exit time probability function associated with the process $X^{(\varepsilon)}$.
\begin{theorem}\label{T:3333}
Under the conditions in Theorem \ref{T:2},
$$
\lim_{\varepsilon\rightarrow 0}\varepsilon\log v_{\varepsilon}(T)=-\inf_{f\in{\cal A}_T^{(1)}}\widetilde{Q}_T(f),
$$
where $\widetilde{Q}_T$ is the rate function defined by (\ref{E:vunzil}).
\end{theorem}

Theorem \ref{T:3333} follows from (\ref{E:exit}) and (\ref{E:chetyre}) with $k=1$.

We will next briefly mention how to use Theorem \ref{T:17} to get a small-noise large deviation style formula for 
the call price function. In the small-noise environment, such a function is defined by 
$C(\varepsilon)=e^{-rT}\mathbb{E}[(S_T^{(\varepsilon)}-K)^{+}]$, $\varepsilon\in(0,1]$,
where $K> 0$ is the strike price and $T$ is the maturity of the option, which are assumed to be fixed. In the previous formula, $u^{+}=
\max(u,0)$, for all $u\in\mathbb{R}$. 

It is well known how to derive a small-noise asymptotic formula for the call price using an 
appropriate large deviation principle (see, e.g., Section 5 in \cite{Ph}; or the proof of Corollary 31 on pages 1130-1133 in \cite{G1};
or the proof of Theorem 5.2 (i) in \cite{G2}). In our case, we will first obtain an asymptotic formula for the price of the binary call option defined by
\begin{equation}
c(\varepsilon)=e^{-rT}\mathbb{P}\left(X_T^{(\varepsilon)}\ge\log K\right),\quad\varepsilon\in(0,1].
\label{E:sdvig}
\end{equation}
Such a formula can be established using the large deviation principle in Theorem \ref{T:17} and the continuity result in Lemma \ref{L:ptu}.
Here we assume that the call option is out-of-the money, that is, $K> s_0e^{rT}$, in order to be able to use Lemma \ref{L:ptu}.
\begin{theorem}\label{T:final}
Suppose the conditions in Theorem \ref{T:17} and Lemma \ref{L:ptu} hold. Suppose also that $K> s_0e^{rT}$. Then
$$
\lim_{\varepsilon\rightarrow 0}\varepsilon\log c(\varepsilon)=-\inf_{x:x\ge\log K-x_0}\widehat{I}_T(x),
$$
where $\widehat{I}_T$ is the rate function defined in Lemma \ref{L:dgh}.
\end{theorem}

\it Proof. \rm The condition $K> s_0e^{rT}$ implies that $rT<\log K-x_0$. Using the continuity of the function $\widehat{I}_T$ on 
the set $\mathbb{R}\backslash rT$ (see Lemma \ref{L:ptu}), we see 
that the set $[\log K-x_0,\infty)$ is a set of continuity for the rate function $\widehat{I}_T$. Now, the equality in Theorem \ref{T:final} can be obtained from (\ref{E:sdvig}) and the large deviation principle in Theorem \ref{T:17}.

This completes the proof of Theorem \ref{T:final}.

We will next explain how to obtain a small-noise large deviation style asymptotic formula for the call price 
$\varepsilon\mapsto C(\varepsilon)$ using Theorem \ref{T:final}. Such a formula can be established by reasoning as in the proof of Corollary 31 in \cite{G1} (see also the proof of part (i) of Theorem 5.2 in \cite{G2}). The proof of the lower large deviation estimate for the call price function uses Theorem \ref{T:final} and Lemma \ref{L:ptu}, while in the proof of the upper estimate, we also assume that the sublinear growth condition holds for the function $\sigma$. This condition is as follows: There exist constants $c_1> 0$ and $c_2> 0$ such that
\begin{equation}
\sigma(t,x)^2\le c_1+c_2x^2,\quad(t,x)\in[0,T]\times\mathbb{R}.
\label{E:lgc}
\end{equation}
\begin{remark}\label{R:mar}
It can be established that the inequality in (\ref{E:lgc}) implies that the discounted asset price process is a martingale, 
and hence $\mathbb{P}$ is a risk-free measure (see similar statements in \cite{FZ,G1}).
\end{remark} 

The upper large deviation estimate for the call price can be obtained by imitating the proof of a similar estimate 
on pages 1131-1133 in \cite{G1}. An important part of the proof of the upper estimate is based on the following statement: There exists 
a constant $\alpha> 0$ such that
\begin{equation}
\mathbb{E}\left[\exp\left\{\alpha\int_0^T\widehat{B}_t^2dt\right\}\right]<\infty.
\label{E:hi}
\end{equation}
The inequality in (\ref{E:hi}) can be established by reasoning as in the proof of Lemma 34 in \cite{G1}. 

Summarizing what was said above, we see that the following theorem is valid.
\begin{theorem}\label{T:tre}
Suppose the conditions in Theorem \ref{T:17} and Lemma \ref{L:ptu} hold. Suppose also that $K> s_0e^{rT}$, and the volatility 
function $\sigma$ satisfies the sublinear growth condition. Then
$$
\lim_{\varepsilon\rightarrow 0}\varepsilon\log C(\varepsilon)=-\inf_{x:x\ge\log K-x_0}\widehat{I}_T(x),
$$
where $\widehat{I}_T$ is the rate function defined in Lemma \ref{L:dgh}.
\end{theorem}
\begin{remark}\label{R:il}
In \cite{G1,G2}, we obtained asymptotic formulas for call prices in the case, where a Gaussian stochastic volatility model is driftless and time-homogeneous, and in addition, the modulus of continuity associated with the volatility process is of power type (see Section 7 of \cite{G1} and Theorem 7.1 (i) in \cite{G2}). The formula in Theorem \ref{T:tre} is valid for more general models, e.g., for the Gaussian logarithmic model with $\beta> 1$, in which the volatility function satisfies the sublinear growth condition. 
\end{remark}
\begin{remark}\label{R:tyw}
The large deviation style formulas obtained in Theorems \ref{T:333}, \ref{T:final}, and \ref{T:tre} provide dominant 
exponential factors in small-noise representations of price functions of various options in super rough Gaussian models.
For example, the formula in (\ref{E:chetyre}) can be rewritten as follows: 
$V_k(\varepsilon)=\exp\{-\frac{L_k}{\varepsilon}\}\exp\{\frac{o(1)}{\varepsilon}\}$ as $\varepsilon\rightarrow 0$.
Therefore, the expression $\exp\{-\frac{L_k}{\varepsilon}\}$ is the dominant exponential factor in a previous 
small-noise representation of the price function of the corresponding binary barrier option. It would be interesting to find explicit 
formulas for the next exponential factors in the above-mentioned representations, or more sharp 
small-noise asymptotic expansions of price functions of options in a super rough Gaussian model. Such expansions were obtained in \cite{BFGHS} for certain fractional Gaussian models. We refer the interested reader to \cite{T} and the references therein for a survey of 
applications of asymptotic expansions in finance.
\end{remark}

\section{Proof of Theorem \ref{T:27}}\label{S:101}
In this section, we provide a detailed proof of Theorem \ref{T:27}. The proof splits into several parts (see Subsections \ref{SS:rv} -- \ref{SS:conti}). 
\subsection{Borel Probability Measures in Banach spaces}\label{SS:rv}
Suppose ${\cal G}$ is a Banach space over the field $\mathbb{R}$ of real numbers. The dual space of ${\cal G}$ is denoted 
by ${\cal G}^{*}$. The space ${\cal G}$ is equipped with the Borel $\sigma$-algebra 
${\cal B}({\cal G})$. The duality between ${\cal G}$ and ${\cal G}^{*}$ will be denoted by $\langle\cdot,\cdot\rangle$. 
For $g\in{\cal G}$ and $\mu\in{\cal G}^{*}$, the symbol $\langle g,\mu\rangle$ stands for the number $\mu(g)$.

Our next goal is to discuss Borel probability measures on the space ${\cal G}$, that is, probability measures on the measurable space
$({\cal G},{\cal B}({\cal G}))$. The first moment of such a measure $\zeta$ is defined by
$
M_1(\zeta)=\int_{{\cal G}}||x||_{\cal G}d\zeta(x)
$
and second moment by
$
M_2(\zeta)=\int_{{\cal G}}||x||_{\cal G}^2d\zeta(x).
$
\begin{definition}\label{D:mvo}
(a)\,Let $\zeta$ be a Borel probability measure on ${\cal G}$ with a finite first moment. If there exists a vector $\bar{m}\in{\cal G}$ such that
$
\langle\bar{m},\mu\rangle=\int_{{\cal G}}\langle x,\mu\rangle d\zeta(x),
$ 
for all $\mu\in{\cal G}^{*}$, then the vector $\bar{m}$ is called the mean vector of $\zeta$. If $\bar{m}=0$, then the measure $\zeta$ is called centered. \\
(b)\,Let $\zeta$ be a Borel probability measure with a finite second moment, and suppose the mean vector $\bar{m}$ exists. 
The covariance of $\zeta$ is the mapping $\mbox{cov}:({\cal G}^{*})^2\mapsto\mathbb{R}$ defined by
$$
\mbox{cov}(\mu_1,\mu_2)=\int_{{\cal G}}\langle x-\bar{m},\mu_1\rangle\langle x-\bar{m},\mu_2\rangle d\zeta(g).
$$
(c)\,Let $\zeta$ be a Borel probability measure with a finite second moment, and suppose the mean vector $\bar{m}$ exists. 
Suppose also that there exists a linear operator $\widehat{K}:{\cal G}^{*}\mapsto{\cal G}$ 
such that for all $\mu_1,\mu_2\in{\cal G}^{*}$,
$\langle\widehat{K}\mu_1,\mu_2\rangle=\mbox{cov}(\mu_1,\mu_2)$.
Then the operator $\widehat{K}$ is called the covariance operator of $\zeta$.
\end{definition}

It is not hard to see that if the mean vector and the covariance operator exist, then they are unique.
The operator $\widehat{K}$ is a bounded operator with
$$
||\widehat{K}||_{{\cal G}^{*}\mapsto{\cal G}}\le\int_{{\cal G}}||x-\bar{m}||_{{\cal G}}^2d\zeta(x)<\infty.
$$
In addition, this operator is symmetric, that is,
$
\langle\widehat{K}\mu_1,\mu_2\rangle=\langle\widehat{K}\mu_2,\mu_1\rangle,
$
for all $\mu_1,\mu_2\in{\cal G}^{*}$.
It is also a non-negative definite operator, i.e., $\langle\widehat{K}\mu,\mu\rangle\ge 0$, 
for all $\mu\in{\cal G}^{*}$. 
\begin{remark}\label{R:above}
It is known that if the space ${\cal G}$ is separable and the conditions in Definition \ref{D:mvo} hold, then the mean vector and the covariance operator exist (see, e.g., \cite{VTC,W} and the references therein).
\end{remark}
\begin{definition}\label{D:rve}
A random vector ${\cal X}$ on $(\Omega,{\cal F},\mathbb{P})$ 
with values in a Banach space ${\cal G}$ is a measurable mapping ${\cal X}:(\Omega,{\cal F})\mapsto({\cal G},{\cal B}({\cal G}))$. 
\end{definition}

A random vector ${\cal X}$ generates a probability measure $\zeta_{{\cal X}}$ as follows: For a set $S\in{\cal B}({\cal G})$, 
$\zeta_{{\cal X}}(S)=\mathbb{P}({\cal X}\in S)$. The measure $\zeta_{{\cal X}}$ is called the distribution of the vector ${\cal X}$.
A random vector ${\cal X}$, for which $M_1\left(\zeta_{{\cal X}}\right)<\infty$ is called integrable, while if 
$M_2\left(\zeta_{{\cal X}}\right)<\infty$, the vector is called square-integrable.

Suppose ${\cal X}$ is a random vector with values in ${\cal G}$, and let $\zeta_{{\cal X}}$ be its distribution. 
The mean vector $\bar{m}$ of the measure $\zeta_{{\cal X}}$ is called the mean vector of ${\cal X}$, while the covariance operator 
$\widehat{K}$ of $\zeta_{{\cal X}}$ is called the covariance operator of ${\cal X}$. Note that the vector $\bar{m}$ and the operator
$\widehat{K}$ associated with the random vector ${\cal X}$ depend on its distribution, and not on the vector itself.

It follows from Remark \ref{R:above} that if the Banach space ${\cal G}$ is separable, then any integrable random vector taking values 
in ${\cal G}$ possesses the mean vector, while 
any square-integrable random vector with values in ${\cal G}$ possesses the covariance operator.

\subsection{Gaussian Measures and Gaussian Vectors}\label{SS:GV}
The following definition introduces Gaussian probability measures and Gaussian random vectors on Banach spaces. 
\begin{definition}\label{D:pf}
Let ${\cal G}$ be a Banach space, and suppose $\zeta$ is a probability measure on the measurable space $({\cal G},{\cal B}({\cal G}))$. 
The measure $\zeta$ is called a Gaussian probability measure if for every $\mu\in{\cal G}^{*}$, the random variable $x\mapsto 
\langle x,\mu\rangle$ is a Gaussian random variable on the measure space $({\cal G},{\cal B}({\cal G}),\zeta)$. The measure $\zeta$ is called non-degenerate if $\int_{{\cal G}}\langle x,\mu\rangle^2d\zeta(x)> 0$ for all $\mu\in{\cal G}^{*}$ except $\mu=0$.
A random vector ${\cal X}$ on a measure space $(\Omega,{\cal F},\mathbb{P})$ with values in ${\cal G}$ 
is called Gaussian if for every $\mu\in{\cal G}^{*}$, the random variable $\langle{\cal X},\mu\rangle$ is normally distributed. 
\end{definition}
It is clear from Definition \ref{D:pf} that a random vector is Gaussian if and only if its distribution is a Gaussian probability measure.

Let ${\cal B}$ be a separable Banach space. It follows from Remark \ref{R:above} that for any Gaussian probability measure $\zeta$ on 
$({\cal G},{\cal B}({\cal G}))$, there exist the mean vector $\bar{m}\in{\cal G}$ and the covariance operator 
$\widehat{K}:{\cal G}^{*}\mapsto{\cal G}$. Moreover, for every $\mu\in{\cal G}^{*}$,
the random variable $M_{\mu}(x)=\langle x,\mu\rangle$ on $({\cal G},{\cal B}({\cal G}),\zeta)$ is normally distributed with mean 
$
m_{\mu}=\langle\bar{m},\mu\rangle
$
and variance 
$
\sigma^2_{\mu}=\langle\widehat{K}\mu,\mu\rangle. 
$
\begin{remark}\label{R:spe}
We refer the reader to \cite{Kuo,Str,VTC} for more information on Borel measures and Gaussian measures on Banach spaces.
\end{remark} 

\subsection{A Special Three-Component Process}\label{SS:spe}
It is known that the dual space 
$\mathbb{C}^{*}$ of $\mathbb{C}$ is the space of all finite signed
Borel measures on $[0,T]$. The norm in the space $\mathbb{C}^{*}$ is defined by
$||\nu||_{\mathbb{C}^{*}}=|\nu|([0,T])$, $\nu\in\mathbb{C}^{*}$, 
where the symbol $|\nu|$ stands for the variation of $\nu$. We will also use the representation $\nu=\nu^{+}-\nu^{-}$,
where $\nu^{+}$ and $\nu^{-}$ are the positive and the negative variations of $\nu$, respectively. Both $\nu^{+}$ and $\nu^{-}$ are finite Borel measures on ${\cal B}([0,T])$. Moreover $|\nu|=\nu^{+}+\nu^{-}$ (see, e.g., \cite{SSh}, Ch. 1, Sect. 7).

Let $\widetilde{\mathbb{C}}$ be a closed subspace of the space $\mathbb{C}$. The dual space $\widetilde{\mathbb{C}}^{*}$
of $\widetilde{\mathbb{C}}$ is the quotient space of $\mathbb{C}^{*}$ by the subspace ${\cal A}$ annihilating $\widetilde{\mathbb{C}}$.
The norm of a coset $S\in\widetilde{\mathbb{C}}^{*}$ is the quotient norm, that is, 
$
||S||_{\widetilde{\mathbb{C}}^{*}}=\inf_{\nu\in S}|\nu|([0,T]),
$
(see, e.g., \cite{Rudin}, Section 4.8). A special example here is the closed subspace $\mathbb{C}_0$ of the space $\mathbb{C}$ 
consisting of all the functions $f$ with $f(0)=0$.
The dual space $(\mathbb{C}_0)^{*}$ of the space $\mathbb{C}_0$ is the quotient space 
of $\mathbb{C}^{*}$ by the annihilator $A$ of $\mathbb{C}_0$. It is not hard to see that the elements of $A$ are the
constant multiples of $\delta_0$, where the symbol $\delta_0$ stands for the $\delta$-measure concentrated at $t=0$.
We will also use the triple direct product $\mathbb{C}^3$ equipped with the norm
$
||(f_1,f_2,f_3)||_{\mathbb{C}^3}=\max\left\{||f_1||_{\mathbb{C}},||f_2||_{\mathbb{C}},||f_3||_{\mathbb{C}}\right\},
$
for all $(f_1,f_2,f_3)\in\mathbb{C}^3$. Its dual space satisfies $(\mathbb{C}^3)^{*}=(\mathbb{C}^{*})^3$.

In the present subsection, we assume that the volatility process $\widehat{B}_t$, $t\in[0,T]$ (see (\ref{E:f11})), 
is a continuous Volterra Gaussian process, for which Assumption A holds. 

Consider the following three-component stochastic process: 
\begin{equation}
{\cal X}_t=(W_t,B_t,\widehat{B}_t),\quad t\in[0,T],
\label{E:ll10}
\end{equation}
with state space $\mathbb{R}^3$, and the associated random vector
${\cal X}:\Omega\mapsto\mathbb{C}^3$ defined by
\begin{equation}
{\cal X}=(W,B,\widehat{B}).
\label{E:ll20}
\end{equation} 
Actually, the vector ${\cal X}$ takes values in a smaller closed subspace of
the space $\mathbb{C}^3$. This will be established below.

We will first show that the random 
vector $\widehat{B}:\Omega\mapsto\mathbb{C}$ takes values in a closed subspace of the space $\mathbb{C}_0$. Consider the mapping 
$\gamma:L^2\mapsto\mathbb{C}$ defined by 
\begin{equation}
\gamma(f)(t)=\int_0^tK(t,s)f(s)ds,\quad t\in[0,T],\quad f\in L^2,
\label{E:yy}
\end{equation} 
and denote 
$
{\cal W}=\overline{\gamma(L^2)},
$
where the closure is taken in the space $\mathbb{C}$. It is not hard to see that ${\cal W}\subset\mathbb{C}_0$.
\begin{lemma}\label{L:takes}
Suppose that Assumption A holds for the kernel of the process $\widehat{B}_t$, $t\in[0,T]$. 
Then $\widehat{B}$ is a Gaussian random vector in the space ${\cal W}$.
\end{lemma}

\it Proof. \rm Denote by ${\cal A}$ the subspace of $\mathbb{C}^{*}$ annihilating ${\cal W}$. Then the dual space ${\cal W}^{*}$ is the quotient space of $\mathbb{C}^{*}$ by ${\cal A}$. To prove that $\widehat{B}$ takes values in ${\cal W}$, it suffices to show that if 
$\nu\in{\cal A}$, then 
\begin{equation}
\int_0^T\widehat{B}_td\nu(t)=0\quad\mbox{on}\quad \Omega.
\label{E:urr}
\end{equation} 
The previous statement follows from the fact that a closed subspace of a Banach space can be separated from a point not belonging to it by a bounded linear functional.

Let $\nu\in{\cal A}$. Then for every $f\in L^2$, 
\begin{equation}
\int_0^Td\nu(t)\int_0^tK(t,s)f(s)ds=0.
\label{E:fr1}
\end{equation}
Next, we see using (\ref{E:mmn}) that
\begin{align*}
\int_0^Td|\nu|(t)\int_0^t|K(t,s)||f(s)|ds&\le||f||_2\int_0^Td|\nu|(t)\left\{\int_0^sK(t,s)^2ds\right\}^{\frac{1}{2}} \\
&\le c||f||_2|\nu|([0,T]),
\end{align*}
where $c> 0$ depends only on $K$. It follows from the previous estimate that Fubini's theorem can be applied to the integral in 
(\ref{E:fr1}). Here we use the following measure spaces:
$([0,T],{\cal B},\nu^{+})$, $([0,T],{\cal B},\nu^{-})$, and $([0,T],{\cal L},l)$, where ${\cal B}$ is the Borel $\sigma$-algebra of $[0,T]$, 
${\cal L}$ is the Lebesgue $\sigma$-algebra of $[0,T]$, and $l$ is the Lebesgue measure on ${\cal L}$.

It follows from Fubini's theorem that
$
\int_0^Tf(s)ds\int_s^TK(t,s)d\nu(t)=0,
$  
for all $f\in L^2$. Therefore, there exists a set $S_{\nu}\in{\cal L}$ such that $l\left(S_{\nu}\right)=T$
and
\begin{equation}
\int_s^TK(t,s)d\nu(t)=0 
\label{E:lpo}
\end{equation}
for all $s\in S_{\nu}$.

Our next goal is to transform the integral $\int_0^T\widehat{B}_td\nu(t)$ using the stochastic Fubini theorem. The function 
$K$ is Borel measurable on $[0,T]^2$. Moreover, using (\ref{E:mmn}), we obtain
$$
\int_0^Td|\nu|(t)\left\{\int_0^tK(t,s)^2ds\right\}^{\frac{1}{2}}\le c|\nu|([0,T])<\infty,
$$
for some $c> 0$. The previous inequality allows us to use the stochastic Fubini theorem 
(see \cite{V} and the references therein).
It follows that
$$
\int_0^T\widehat{B}_td\nu(t)=\int_0^Td\nu(t)\int_0^tK(t,s)dB_s=\int_0^TdB_s\int_s^TK(t,s)d\nu(t).
$$
Now, (\ref{E:lpo}) implies that $\int_0^T\widehat{B}_td\nu(t)=0$ on $\Omega$.

This establishes (\ref{E:urr}) and completes the proof of Lemma \ref{L:takes}.
\begin{remark}\label{R:embed}
The mapping $\gamma$ in (\ref{E:yy}) is continuous, but not necessarily an embedding. 
\end{remark}

We will next return to the random process ${\cal X}_t$, $t\in[0,T]$, with state space $\mathbb{C}^3$ and the associated random vector 
${\cal X}$ in $\mathbb{C}^3$, defined in (\ref{E:ll10}) and (\ref{E:ll20}), respectively.
\begin{theorem}\label{T:hg}
Let $W_t$, $t\in[0,T]$, be a standard Brownian motion on $(\Omega, {\cal F},\mathbb{P})$ independent of $B_t$, $t\in[0,T]$, 
and suppose Assumption A holds. 
Then ${\cal X}=(W,B,\widehat{B})$ is a centered Gaussian random vector in
the space $\mathbb{C}^3$. 
\end{theorem}
\begin{remark}\label{R:remt}
Using linear combinations of $\delta$-measures, we see that Theorem \ref{T:hg} implies the following statement.
For all systems $0\le t_1<\cdots< t_n\le T$, $0\le s_1<\cdots< s_m\le T$, and $0\le r_1<\cdots< r_p\le T$, the random variable
$(W_{t_1},\cdots,W_{t_n},B_{s_1},\cdots,B_{s_m},\widehat{B}_{r_1},\cdots,\widehat{B}_{r_p})$ is multivariate Gaussian.
\end{remark}

\it Proof of Theorem \ref{T:hg}. \rm The random vector ${\cal X}$ takes values in the space $\mathbb{C}^3$. 
We will next show that the vector ${\cal X}$ is Gaussian. It suffices to prove that for every
$(\lambda_0,\lambda_1,\lambda_2)\in(\mathbb{C}^{*})^3$, the random variable
$
\langle W,\lambda_0\rangle+\langle B,\lambda_1\rangle+\langle\widehat{B},\lambda_2\rangle
$
is Gaussian. In the previous expression, the symbol $\langle\cdot,\cdot\rangle$ stands for the duality between $\mathbb{C}$ 
and $\mathbb{C}^{*}$. 
It is not hard to see that
$
\langle W,\lambda_0\rangle=\int_0^TW_td\lambda_0(t)
$
is a Gaussian random variable. Since $W$ and $B$ are independent, we only need to show that the random variable
$
F:=\langle B,\lambda_1\rangle+\langle\widehat{B},\lambda_2\rangle
$
is Gaussian. We have 
$
F=\int_0^Td\lambda_1(t)\int_0^tdB_s+\int_0^Td\lambda_2(t)\int_0^tK(t,s)dB_s.
$
Applying the stochastic Fubini theorem to the integrals on the right-hand side of the previous equality 
(see the reference in the proof of 
Lemma \ref{L:takes}), we obtain
\begin{equation}
F=\int_0^T\left(\int_s^Td\lambda_1(t)+\int_s^TK(t,s)d\lambda_2(t)\right)dB_s.
\label{E:inte}
\end{equation}
Here it is important to recall that $K$ is a Borel measurable function of two variables. Set
\begin{equation}
A(s;\lambda_1,\lambda_2):=\int_s^Td\lambda_1(t)+\int_s^TK(t,s)d\lambda_2(t),\quad s\in[0,T].
\label{E:gt}
\end{equation}

Our next goal is to analyze the function on the right-hand side of (\ref{E:gt}). 
It is clear that the function 
$s\mapsto\int_s^Td\lambda_1(t)=\lambda_1[s,T]$, $s\in[0,T]$,
is Borel measurable. Moreover, Tonelli's theorem and (\ref{E:mmn}) imply that
\begin{align*}
&\int_0^Tds\int_s^T|K(t,s)|d|\lambda_2|(t)=\int_0^Td|\lambda_2|(t)\int_0^t|K(t,s)|ds \\
&\le||\lambda_2||\sup_{t\in[0,T]}\int_0^t|K(t,s)|ds\le||\lambda_2||T^{\frac{1}{2}}
\sup_{t\in[0,T]}\left\{\int_0^tK(t,s)^2ds\right\}^{\frac{1}{2}}<\infty.
\end{align*}
Now, the Fubini theorem for Borel measures (see a general Fubini-Tonelli theorem for Radon Products in \cite{Fol}, (7.27)) 
implies that the function
$s\mapsto\int_s^TK(t,s)d\lambda_2(t)$, $s\in[0,T]$,
is Borel measurable. Summarizing what was said above, we see that the function $A$ is Borel measurable.

It is also true that the function $A$ is square-integrable over $[0,T]$. Indeed,
$$
|A(s;\lambda_1,\lambda_2)|\le||\lambda_1||+\int_s^T|K(t,s)|d|\lambda_2|(t).
$$
and hence
\begin{align*}
&\int_0^TA(s;\lambda_1,\lambda_2)^2ds\le 2T||\lambda_1||^2+2||\lambda_2||\int_0^Tds\int_s^TK(t,s)^2d|\lambda_2|(t) 
\nonumber \\
&=2T||\lambda_1||^2+2||\lambda_2||\int_0^Td|\lambda_2|(t)\int_0^tK(t,s)^2ds  \\
&\le 2T||\lambda_1||^2+2||\lambda_2||^2\sup_{t\in[0,T]}\int_0^tK(t,s)^2ds
\le 2T||\lambda_1||^2+2C||\lambda_2||^2<\infty.
\end{align*}
In the previous estimates, we used (\ref{E:mmn}). Therefore, the function $A$ is square-integrable over $[0,T]$.
It follows that the stochastic integral on the right-hand side of 
(\ref{E:inte}) exists. It is known that Wiener integrals with square-integrable integrands are Gaussian random variables. 
Therefore, the random variable $F$ is Gaussian.

This completes the proof of Theorem \ref{T:hg}.
 
Theorem \ref{T:hg} states that the random vector ${\cal X}$ is a Gaussian vector in the space $\mathbb{C}^{3}$.
We will next find a smaller space in which the Gaussian vector ${\cal X}$ takes its values.

The next straightforward lemma will be useful in the sequel.
\begin{lemma}\label{L:simm}
Suppose that the conditions in Theorem \ref{T:hg} hold. Suppose also that the range of the mapping ${\cal X}:\Omega\mapsto\mathbb{C}^3$ is a subset of a closed subspace 
$S$ of the space $\mathbb{C}^3$.  Then ${\cal X}$ is a Gaussian random vector in $S$.
\end{lemma}

\it Proof. \rm Let $\widehat{\lambda}\in S^{*}$. Then $\widehat{\lambda}$ is an element of the quotient space of $\mathbb{C}^3$ 
by the annihilator of $S$. It follows that
$\widehat{\lambda}({\cal X})=\lambda({\cal X})$, for every $\lambda\in\widehat{\lambda}$. Therefore the random variable 
$\widehat{\lambda}({\cal X})$ is Gaussian, by Theorem \ref{T:hg}. 

This establishes Lemma \ref{L:simm}.

It is well-known that the Cameron-Martin space $\mathbb{H}_0^1$ is continuously embedded into the space $\mathbb{C}_0$. Moreover, the range of this embedding is dense 
in $\mathbb{C}_0$. Set $H=L^2\times L^2$. Then $H$ is a separable Hilbert space equipped with the norm defined as follows: 
$||h||_H=\sqrt{||h_0||_2^2+||h_1||_2^2}$, $h=(h_0,h_1)\in H$.
Let $j:H\mapsto\mathbb{C}\times\mathbb{C}\times{\cal W}$ be the mapping defined on the space $H$ by 
$j(h_0,h_1)=(g_0,g_1,g_2),$
where for all $t\in[0,T]$, $g_0(t)=\int_0^th_0(s)ds$, $g_1(t)=\int_0^th_1(s)ds$, and 
$g_2(t)=\int_0^tK(t,s)h_1(s)ds$. 
It is easy to see that $j$ is a continuous linear mapping. It is also an injection, since the mappings $h_0\mapsto g_0$ and 
$h_1\mapsto g_1$ are injections. Using the previous formulas, we see that the range $j(H)$ of the mapping $j$ consists of all the triples 
$(g_0,g_1,g_2)$ with $g_0\in H_0^1$, $g_1\in H_0^1$, and 
$g_2(t)=\int_0^tK(t,s)\dot{g}_1(s)ds$, $t\in[0,T]$.

Define a subspace of the space $\mathbb{C}^3$ by
$\widetilde{{\cal G}}=\overline{j(H)}$,
where the closure is taken in the space $(\mathbb{C})^3$. Then $\widetilde{{\cal G}}$ is a separable Banach space, and moreover 
$\widetilde{{\cal G}}\subset\mathbb{C}_0\times\mathbb{C}_0\times{\cal W}$.
The dual space $\widetilde{{\cal G}}^{*}$ of $\widetilde{{\cal G}}$ is the quotient space of the space 
$(\mathbb{C})^3$ by its subspace ${\cal A}$ that
annihilates $j(H)$ ($j(H)$ is dense in $\widetilde{{\cal G}}$). The annihilation condition is the following: For $\alpha=(\alpha_0,
\alpha_1,\alpha_2)\in{\cal A}$,
\begin{align}
\int_0^Td\alpha_0(t)\int_0^th_0(s)ds&+\int_0^Td\alpha_1(t)\int_0^th_1(s)ds
+\int_0^Td\alpha_2(t)\int_0^tK(t,s)h_1(s)ds=0,
\label{E:alpha}
\end{align}
for all functions $h_0\in L^2$ and $h_1\in L^2$. Next, Fubini's theorem implies that an equivalent form of the equality in 
(\ref{E:alpha}) is as follows:
\begin{align}
&\int_0^Th_0(s)ds\int_s^Td\alpha_0(t)+\int_0^Th_1(s)ds\left[\int_s^Td\alpha_1(t)
+\int_s^TK(t,s)d\alpha_2(t)\right]=0,
\label{E:alphas}
\end{align}
for all functions $h_0\in L^2$ and $h_1\in L^2$. Now, plugging $h_1=0$ and then $h_0=0$ into (\ref{E:alphas}), we see that the
equality in (\ref{E:alphas}) splits into the following two equalities:
\begin{equation}
\int_s^Td\alpha_0(t)=0
\label{E:88}
\end{equation}
and
\begin{equation}
\int_s^Td\alpha_1(t)+\int_s^TK(t,s)d\alpha_2(t)=0,
\label{E:99}
\end{equation}
for almost all $s\in[0,T]$ with respect to the Lebesgue measure on $[0,T]$.

The proof of the next assertion resembles that of Lemma \ref{L:takes}.
\begin{theorem}\label{T:wider}
Let $W_t$, $t\in[0,T]$, be a standard Brownian motion on $(\Omega,{\cal F},\mathbb{P})$ independent of $B_t$, $t\in[0,T]$, and 
suppose that Assumption A holds.
Then ${\cal X}=(W,B,\widehat{B})$ is a centered Gaussian random vector in the space $\widetilde{{\cal G}}$. 
\end{theorem}

\it Proof. \rm We will first prove that for every 
$
\alpha=(\alpha_0,\alpha_1,\alpha_2)\in{\cal A},
$
we have
\begin{align}
&\int_0^Td\alpha_0(t)\int_0^tdW_s+\int_0^Td\alpha_1(t)\int_0^tdB_s
+\int_0^Td\alpha_2(t)\int_0^tK(t,s)dB_s=0.
\label{E:alphak}
\end{align}
Using the stochastic Fubini theorem as in Lemma \ref{L:takes}, we see that the equality in (\ref{E:alphak}) is equivalent to 
the equality
\begin{align}
&\int_0^TdW_s\int_s^Td\alpha_0(t)+\int_0^TdB_s\left[\int_s^Td\alpha_1(t)+\int_s^TK(t,s)d\alpha_2(t)\right]=0.
\label{E:alpka}
\end{align} 
It is clear that (\ref{E:alpka}) follows from (\ref{E:88}) and (\ref{E:99}). This proves that the range of the mapping 
${\cal X}:\Omega\mapsto\mathbb{C}^3$ is contained in the space $\widetilde{{\cal G}}$. We have already 
established that ${\cal X}$ is a Gaussian random vector in the space $\mathbb{C}^3$. Now,
Theorem \ref{T:wider} follows from Lemma \ref{L:simm}.

\subsection{Covariance Function and Covariance Operator}\label{SS:cco}
In the present subsection, we compute the covariance function and the covariance operator of the Gaussian random vector 
${\cal X}$ considered in the previous subsection. 
\begin{theorem}\label{T:ccd}
Suppose the process $\widehat{B}_t$, $t\in[0,T]$, in (\ref{E:yy}) is a Gaussian Volterra process. Suppose also that the restrictions in Theorem \ref{T:hg} 
hold, and consider the random vector ${\cal X}=(W,B,\widehat{B})$ as a vector taking values in the space $\mathbb{C}^3$. Then, for all 
$\mu_1=(\lambda_0,\lambda_1,\lambda_2)\in(\mathbb{C}^{*})^3$ and $\mu_2=(\nu_0,\nu_1,\nu_2)\in(\mathbb{C}^{*})^3$,
the covariance of ${\cal X}$ is given by
\begin{align}
\mbox{cov}(\mu_1,\mu_2)&=\int_0^Tf_0(u;\lambda_0,\lambda_1,\lambda_2)d\nu_0(u)
+\int_0^Tf_1(u;\lambda_0,\lambda_1,\lambda_2)d\nu_1(u) \nonumber \\
&\quad+\int_0^Tf_2(u;\lambda_0,\lambda_1,\lambda_2)d\nu_2(u),
\label{E:covr}
\end{align}
where the functions $f_0$, $f_1$, and $f_2$ are defined on $[0,T]$ by
\begin{align*}
&f_0(u;\lambda_0,\lambda_1,\lambda_2)=\int_0^T(t\wedge u)d\lambda_0(t), \\
&f_1(u;\lambda_0,\lambda_1,\lambda_2)=\int_0^T(t\wedge u)d\lambda_1(t)
+\int_0^Td\lambda_2(t)\int_0^{t\wedge u}K(t,s)ds, \\
&f_2(u;\lambda_0,\lambda_1,\lambda_2)=\int_0^u\left(\int_s^Td\lambda_1(t)+\int_s^TK(t,s)d\lambda_2(t)\right)K(u,s)ds.
\end{align*}
In addition, the covariance operator is given by
$
\widehat{K}(\mu)=(f_0,f_1,f_2),
$
for all $\mu=(\lambda_0,\lambda_1,\lambda_2)\in(\mathbb{C}^{*})^3$.
\end{theorem}

\it Proof. \rm  Let $\mu_1=(\lambda_0,\lambda_1,\lambda_2)\in(\mathbb{C}^{*})^3$ and $\mu_2=(\nu_0,\nu_1,\nu_2)\in(\mathbb{C}^{*})^3$. 
Then we have
\begin{align*}
&\mbox{cov}(\mu_1,\mu_2)=\mathbb{E}\left[\langle{\cal X},\mu_1\rangle\langle{\cal X},\mu_2\rangle\right] 
=\mathbb{E}\left[\widetilde{A}(\lambda_0,\lambda_1,\lambda_2)\widetilde{A}(\nu_0,\nu_1,\nu_2)\right],
\end{align*}
where for $(\xi_0,\xi_1,\xi_2)\in(\mathbb{C}^{*})^3$, 
\begin{align*}
&\widetilde{A}(\xi_0,\xi_1,\xi_2)=\int_0^Td\xi_0(t)\int_0^tdW_s+\int_0^Td\xi_1(t)\int_0^tdB_s 
+\int_0^Td\xi_2(t)\int_0^tK(t,s)dB_s.
\end{align*}
Next, using the stochastic Fubini theorem and the independence of $W$ and $B$, we obtain
\begin{align}
\mbox{cov}(\mu_1,\mu_2)&=\int_0^Tds\int_s^Td\lambda_0(t)\int_s^Td\nu_0(u) 
+\int_0^Tds\left(\int_s^Td\lambda_1(t)+
\int_s^TK(t,s)d\lambda_2(t)\right) \nonumber \\
&\quad\times\left(\int_s^Td\nu_1(u)+
\int_s^TK(u,s)d\nu_2(u)\right).
\label{E:E1}
\end{align}
It follows that
\begin{align*}
&\mbox{cov}(\mu_1,\mu_2)=\int_0^Tds\int_s^Td\lambda_0(t)\int_s^Td\nu_0(u) 
+\int_0^Tds\int_s^Td\lambda_1(t)\int_s^Td\nu_1(u) \\
&+\int_0^Tds\int_s^Td\nu_1(u)\int_s^TK(t,s)d\lambda_2(t) 
+\int_0^Tds\int_s^Td\lambda_1(t)\int_s^TK(u,s)d\nu_2(u) \\
&+\int_0^Tds\int_s^TK(t,s)d\lambda_2(t)\int_s^TK(u,s)d\nu_2(u).
\end{align*}

Finally, applying Fubini's theorem, we get
\begin{align*}
&\mbox{cov}(\mu_1,\mu_2)=\int_0^Td\nu_0(u)\int_0^T(t\wedge u)d\lambda_0(t) \\
&+\int_0^Td\nu_1(u)\left[\int_0^T(t\wedge u)d\lambda_1(t)
+\int_0^Td\lambda_2(t)\int_0^{t\wedge u}K(t,s)ds\right] \\
&+\int_0^Td\nu_2(u)\left[\int_0^uK(u,s)ds\int_s^Td\lambda_1(t)
+\int_0^uK(u,s)ds\int_s^TK(t,s)d\lambda_2(t)\right].
\end{align*}
This establishes the formulas for the covariance function and the covariance operator in Theorem \ref{T:ccd}.

The proof of Theorem \ref{T:ccd} is thus completed.
\begin{remark}\label{R:rr}
The functions $f_0$, $f_1$, and $f_2$ appearing in Theorem \ref{T:ccd} are continuous. This can be shown using (\ref{E:mmn}). 
\end{remark}

In the next statement, we use the fact that the vector ${\cal X}$ takes values in the space
$\widetilde{{\cal G}}$ defined in Subsection \ref{SS:spe}. Recall that the dual space $\widetilde{{\cal G}}^{*}$ of the space 
${\cal G}$ is the quotient space of $\mathbb{C}^3$ by the annihilator ${\cal A}$ of ${\cal G}$. Let 
$\widehat{\mu}_1\in\widetilde{{\cal G}}^{*}$ and $\widehat{\mu}_2\in\widetilde{{\cal G}}^{*}$, and denote by 
${\rm Cov}(\widehat{\mu}_1,\widehat{\mu}_2)$ the value of the covariance function. Every element $\widehat{\mu}\in\widetilde{{\cal G}}^{*}$ 
is a coset in $\mathbb{C}^3$, and the notation $\mu\in\widehat{\mu}$ will mean that $\mu\in\mathbb{C}^3$ 
is an element of $\widehat{\mu}$.
\begin{theorem}\label{T:cc}
Suppose the process $\widehat{B}$ in (\ref{E:yy}) is a Gaussian Volterra process. Suppose also that the assumptions in Theorem \ref{T:hg} 
hold, and consider ${\cal X}=(W,B,\widehat{B})$ as a Gaussian random vector in the space $\widetilde{{\cal G}}$. Then, for all 
$\widehat{\mu}_1\in\widetilde{{\cal G}}^{*}$ and $\widehat{\mu}_2\in\widetilde{{\cal G}}^{*}$, the following formula holds:
\begin{equation}
\mbox{\rm Cov}(\widehat{\mu}_1,\widehat{\mu}_2)=\mbox{\rm cov}(\mu_1,\mu_2),
\quad\mu_1\in\widehat{\mu}_1,\quad\mu_2\in\widehat{\mu}_2.
\label{E:opo}
\end{equation}
The covariance $\mbox{\rm cov}$ appearing in formula (\ref{E:opo})
is described in Theorem \ref{T:ccd}, and the value of the expression on the right-hand side
of (\ref{E:opo}) is the same for all $\mu_1\in\widehat{\mu}_1$ and $\mu_2\in\widehat{\mu}_2$.
In addition, the covariance operator $\widehat{K}:\widetilde{{\cal G}}^{*}\mapsto\widetilde{{\cal G}}$ is given by
\begin{equation}
\widehat{K}(\widehat{\mu})=(f_0,f_1,f_2),\quad\widehat{\mu}\in\widetilde{{\cal G}}^{*},
\label{E:tttu}
\end{equation}
where the functions $f_0$, $f_1$, and $f_2$ are defined in the formulation of Theorem \ref{T:ccd}. These functions
satisfy the following condition: If $(\lambda_0,\lambda_1,\lambda_2)\in\widehat{\mu}$ and 
$(\nu_0,\nu_1,\nu_2)\in\widehat{\mu}$, then 
\begin{equation}
f_i(u;\lambda_0,\lambda_1,\lambda_2)=f_i(u;\nu_0,\nu_1,\nu_2),
\label{E:ii}
\end{equation}
for all $u\in[0,T]$ and $0\le i\le 2$.
\end{theorem}

\it Proof. \rm  Let $\widehat{\mu}_1\in\widetilde{{\cal G}}^{*}$ and $\widehat{\mu}_2\in\widetilde{{\cal G}}^{*}$.
Then
\begin{equation}
\mbox{\rm Cov}(\widehat{\mu}_1,\widehat{\mu}_2)=\mathbb{E}\left[\widehat{\mu}_1({\cal X})\widehat{\mu}_2({\cal X})\rangle\right].
\label{E:deef}
\end{equation}
We also have
\begin{equation}
\mbox{\rm cov}(\mu_1,\mu_2)=\mathbb{E}\left[\mu_1({\cal X})\mu_2({\cal X})\right],
\label{E:def}
\end{equation}
for all $\mu_1,\mu_2\in\mathbb{C}^{*}$.

It follows from Theorem \ref{T:wider} that the random vector ${\cal X}=(W,B,\widehat{B})$ takes values in the space 
$\widetilde{{\cal G}}$. Since ${\cal A}$ annihilates $\widetilde{{\cal G}}$, the equalities in (\ref{E:deef}) and (\ref{E:def}) 
imply that (\ref{E:opo}) holds. Moreover,
it follows from (\ref{E:def}) that for every $\mu_1=(\lambda_0,\lambda_1,\lambda_2)\in(\mathbb{C}^3)^{*}$, 
$\mu_2=(\nu_0,\nu_1,\nu_2)\in(\mathbb{C}^3)^{*}$, and $\alpha=(\alpha_0,\alpha_1,\alpha_2)\in{\cal A}$,
\begin{equation}
\mbox{\rm cov}(\mu_1,\mu_2)=\mbox{\rm cov}(\mu_1+\alpha,\mu_2).
\label{E:label}
\end{equation}
Now, the equality in (\ref{E:ii}) follows from (\ref{E:covr}). Indeed, by taking into account (\ref{E:label}) and 
plugging $\nu_1=\nu_2=0$ into (\ref{E:covr}), we see that the equality
$$
\int_0^Tf_0(u;\lambda_0+\alpha_0,\lambda_1+\alpha_1,\lambda_2+\alpha_2)d\nu_0(u)
=\int_0^Tf_0(u;\lambda_0,\lambda_1,\lambda_2)d\nu_0(u)
$$
holds for all $\nu_0\in\mathbb{C}^{*}$. It follows that 
$
f_0(u;\lambda_0+\alpha_0,\lambda_1+\alpha_1,\lambda_2+\alpha_2)=f_0(u;\lambda_0,\lambda_1,\lambda_2).
$
This establishes (\ref{E:ii}) for the function $f_0$. The proof of (\ref{E:ii}) for the functions
$f_1$ and $f_2$ is similar.

Finally, the formula for the covariance operator in (\ref{E:tttu}) follows from Theorem \ref{T:ccd}.

The proof of Theorem \ref{T:cc} is thus completed.

\subsection{An Abstract Wiener Space Associated with the Random Vector ${\cal X}$}\label{SS:ert}
The notion of an abstract Wiener space goes back to Gross (see \cite{Gross}). We use the definition in \cite{BBAK}.
\begin{definition}\label{D:aws}
An abstract Wiener space is a quadruple $({\cal G},H,j,\zeta)$, where 
\begin{enumerate}
\item ${\cal G}$ is a separable Banach space.
\item $H$ is a separable Hilbert space.
\item $j:H\mapsto{\cal G}$ is a continuous embedding (linear injection) such that 
$j(H)$ is dense in ${\cal G}$.
\item $\zeta$ is a probability measure on $({\cal G},{\cal B}({\cal G}))$ such that
for every $\mu\in{\cal G}^{*}$,
\begin{equation}
\int_{{\cal G}}\exp\{i\langle x,\mu\rangle\}d\zeta(x)=\exp\left\{-\frac{1}{2}||j^{*}(\mu)||_H^2\right\}.
\label{E:rew}
\end{equation}
\end{enumerate}
In (\ref{E:rew}), the symbol $j^{*}$ stands for the adjoint transformation $j^{*}:{\cal G}^{*}\mapsto H^{*}=H$.
\end{definition}

The next statement is standard (see, e.g., \cite{B}).
\begin{lemma}\label{L:norm}
The characteristic function of a normally distributed random variable $X$ with mean $m$ and variance $\sigma^2$ is given by the following formula: $\mathbb{E}\left[e^{i\xi X}\right]=\exp\left\{im\xi-\frac{\sigma^2\xi^2}{2}\right\}$, $\xi\in\mathbb{R}$.
Conversely, if the characteristic function of a random variable $X$ satisfies the previous equality, then $X$ is 
$N(m,\sigma^2)$-distributed.
\end{lemma}
\begin{remark}\label{R:rth}
By replacing $\mu$ in (\ref{E:rew}) by $\xi\mu$ with $\xi\in\mathbb{R}$, and using Lemma \ref{L:norm}, we see that the measure $\zeta$ in item 4 of Definition
\ref{D:aws} is a centered Gaussian measure on $({\cal G},{\cal B}({\cal G}))$.
\end{remark}

In the proof of Theorem \ref{T:2}, the following large deviation principle for abstract Wiener spaces will be used (see, e.g., Theorem 3.4.12 in \cite{DS}, Theorem 2.3 in \cite{BBAK}, Theorem 8.4.1 in \cite{Str}).
\begin{theorem}\label{T:lpd1}
Let $({\cal G},H,j,\zeta)$ be an abstract Wiener space. Then for every Borel subset $A\subset{\cal G}$,
\begin{align*}
-\inf_{x\in A^{\circ}}\lambda_{\zeta}(x)\le\liminf_{\varepsilon\rightarrow 0}\varepsilon^2\log\zeta\left(\frac{1}{\varepsilon}A\right)
\le\limsup_{\varepsilon\rightarrow 0}\varepsilon^2\log\zeta\left(\frac{1}{\varepsilon}A\right)\le-\inf_{x\in\bar{A}}\lambda_{\zeta}(x),
\end{align*}
where the rate function $\lambda_{\zeta}:{\cal G}\mapsto[0,\infty]$ is defined by
$$
\lambda_{\zeta}(x)=\begin{cases}
\frac{1}{2}||j^{-1}x||_H, & x\in j(H) \\
\infty, & x\notin j(H).
\end{cases}
$$
\end{theorem}

We will next construct a special abstract Wiener space associated with the random vector ${\cal X}$. 
Recall that in Subsection \ref{SS:spe} we used the notation $H=L^2\times L^2$, and defined the mapping $j:H\mapsto\mathbb{C}\times\mathbb{C}\times{\cal W}$ by 
$j(h_0,h_1)=(g_0,g_1,g_2)$,
where $g_0(t)=\int_0^th_0(s)ds$, $g_1(t)=\int_0^th_1(s)ds$, and 
$g_2(t)=\int_0^tK(t,s)h_1(s)ds$, $t\in[0,T]$. We also set 
$\widetilde{{\cal G}}=\overline{j(H)}$,
where the closure is taken in the space $(\mathbb{C})^3$.

Define the probability measure $\zeta$ on 
$(\widetilde{{\cal G}},{\cal B}(\widetilde{{\cal G}}))$ by 
$\zeta(A)=\mathbb{P}((W,B,\widehat{B})\in A)$, $A\in{\cal B}(\widetilde{{\cal G}})$. 
By Theorem \ref{T:wider}, $(W,B,\widehat{B})$ is a centered Gaussian vector taking values in the space $\widetilde{{\cal G}}$. 
Therefore, the probability measure $\zeta$ defined above is a Gaussian measure
satisfying
$$
\int_{\widetilde{{\cal G}}}\exp\{i\langle x,\widehat{\mu}\rangle\}d\zeta(x)=\exp\left\{-\frac{1}{2}\langle\widehat{K}\widehat{\mu},
\widehat{\mu}\rangle\right\},\quad\widehat{\mu}\in\widetilde{{\cal G}}^{*},
$$
where the covariance operator $\widehat{K}$ is described in Theorem \ref{T:cc}.
\begin{theorem}\label{T:aVs}
The quadruple $(\widetilde{{\cal G}},H,j,\zeta)$ associated with the random vector ${\cal X}=(W,B,\widehat{B})$ is an abstract Wiener space.
\end{theorem}

\it Proof. \rm We only need to establish that for all $\widehat{\mu}\in\widetilde{{\cal G}}^{*}$,
\begin{equation}
\langle\widehat{K}\widehat{\mu},\widehat{\mu}\rangle=||j^{*}(\widehat{\mu})||_H^2.
\label{E:trr}
\end{equation}
The functional $\widehat{\mu}\in\widetilde{{\cal G}}^{*}$ is an element of the quotient space of $(\mathbb{C}^{*})^3$ by ${\cal A}$.
Let $\mu=(\lambda_0,\lambda_1,\lambda_2)\in\widehat{\mu}$. Using Theorem \ref{T:ccd} and (\ref{E:E1}), we obtain
\begin{align}
&\langle\widehat{K}\widehat{\mu},\widehat{\mu}\rangle=\mbox{\rm cov}(\mu,\mu)
=\int_0^Tds\left[\int_s^Td\lambda_0(t)\right]^2 \nonumber \\
&\quad+\int_0^Tds\left[\int_s^Td\lambda_1(t)
+\int_s^TK(t,s)d\lambda_2(t)\right]^2.
\label{E:cova1}
\end{align}

Our next goal is to compute $j^{*}(\widehat{\mu})$. For all $g=(h_0,h_1)\in H$, we have
\begin{align*}
&j^{*}(\widehat{\mu})(g)=\widehat{\mu}(j(g))=\mu(j(g)) \\
&=\int_0^Td\lambda_0(t)\int_0^th_0(s)ds+\int_0^Td\lambda_1(t)\int_0^t
h_1(s)ds+\int_0^Td\lambda_2(t)\int_0^tK(t,s)h_1(s)ds \\
&=\int_0^Th_0(s)ds\int_s^Td\lambda_0(t)+\int_0^Th_1(s)ds\left(\int_s^Td\lambda_1(t)+\int_s^TK(t,s)d\lambda_2(t)\right),
\end{align*}
for all $\mu=(\lambda_0,\lambda_1,\lambda_2)\in\widehat{\mu}$. It follows that
$$
j^{*}(\widehat{\mu})=\left(\int_{\cdot}^Td\lambda_0(t),\int_{\cdot}^Td\lambda_1(t)
+\int_{\cdot}^TK(t,s)d\lambda_2(t)\right).
$$
Therefore,
\begin{align}
&||j^{*}(\widehat{\mu})||_H=\int_0^Tds\left[\int_s^Td\lambda_0(t)\right]^2 
+\int_0^Tds\left[\int_s^Td\lambda_1(t)
+\int_s^TK(t,s)d\lambda_2(t)\right]^2.
\label{E:cova2}
\end{align}
Now, (\ref{E:trr}) follows from (\ref{E:cova1}) and (\ref{E:cova2}).

This completes the proof of Theorem \ref{T:aVs}.

Our next goal is to prove a large deviation principle for the special abstract Wiener space constructed in 
Theorem \ref{T:aVs}.
\begin{theorem}\label{T:ladp}
Suppose Assumption A holds. Then the following estimates are valid for the random vector 
${\cal X}=(W,B,\widehat{B})$:
For every Borel subset $E$ of ${\cal G}$,
\begin{align*}
-\inf_{(g_0,g_1,g_2)\in E^{\circ}}\Lambda(g_0,g_1,g_2)&\le\liminf_{\varepsilon\rightarrow 0}\varepsilon
\log\mathbb{P}\left(\sqrt{\varepsilon}{\cal X}\in E\right)
\le\limsup_{\varepsilon\rightarrow 0}\varepsilon\log\mathbb{P}\left(\sqrt{\varepsilon}{\cal X}\in E\right) \\
\le-\inf_{(g_0,g_1,g_2)\in\overline{E}}\Lambda(g_0,g_1,g_2),
\end{align*}
where the good rate function $\Lambda:{\cal G}\mapsto[0,\infty]$ is defined by
$$
\Lambda(g_0,g_1,g_2)=\frac{1}{2}\int_0^T\dot{g}_0(s)^2ds+\frac{1}{2}\int_0^T\dot{g}_1(s)^2ds,
$$
for all $(g_0,g_1,g_2)\in{(\mathbb{H}^1_0)^3}$ with $g_2=\widehat{g}_1$, while $\Lambda(g_0,g_1,g_2)=\infty$
otherwise.
\end{theorem}

\it Proof. \rm  For the abstract Wiener space in Theorem \ref{T:ladp}, 
we have 
$$
j(H)=\{(g_0,g_1,g_2)\in(\mathbb{H}_0^1)^3:g_2=\widehat{g}_1\}.
$$ 
Moreover if $(g_0,g_1,g_2)\in j(H)$, then $j^{-1}(g_0,g_1,g_2)=(h_0,h_1)\in L^2\times L^2$. In the previous equality,
$h_0(t)=\dot{g}_0(t)$ and $h_1(t)=\dot{g}_1(t)$, for all $t\in[0,T]$.
Therefore,
$$
||j^{-1}(g_0,g_1,g_2)||_H^2=\frac{1}{2}\int_0^T\dot{g}_0(s)^2ds+\frac{1}{2}\int_0^T\dot{g}_1(s)^2ds.
$$

Now, it is clear that Theorem \ref{T:ladp} follows from Theorem \ref{T:lpd1}.
\begin{remark}\label{R:foll}
Theorem \ref{T:ladp} and the contraction principle imply that the process $\varepsilon\mapsto\sqrt{\varepsilon}\widehat{B}$ 
satisfies a sample path large deviation principle
with speed $\varepsilon^{-1}$ and good rate function $J$ given on the space $\mathbb{C}_0[0,1]$ by
$
J(g)=\frac{1}{2}\inf_{\{f\in\mathbb{H}_0^1:\widehat{f}=g\}}\int_0^T\dot{f}(t)^2dt,
$
if the equation $\widehat{f}=g$ is solvable for $f$, and $J(g)=\infty$, otherwise.
\end{remark}

The next corollary is a key ingredient in the proof of Theorem \ref{T:27} (see Subsection \ref{SS:conti}). 
\begin{corollary}\label{C:corol}
Let $\widehat{B}_t$, $t\in[0,T]$, be a Gaussian Volterra volatility process such that Assumption A holds for its kernel $K$. 
Then the following equality is valid:
\begin{equation}
\limsup_{m\rightarrow\infty}\limsup_{\varepsilon\rightarrow 0}\varepsilon
\log\mathbb{P}(\sup_{t,s\in[0,T]:|t-s|<\frac{T}{m}}|\widehat{B}_t-\widehat{B}_s|\ge\varepsilon^{-\frac{1}{2}}y)
=-\infty,
\label{E:cvf}
\end{equation}
for all $y> 0$, where $m$ denotes a positive integer, and $\varepsilon> 0$ is a small parameter. 
\end{corollary}

\it Proof. \rm Corollary \ref{C:corol} will be derived from the LDP in Remark \ref{R:foll}. Let $y> 0$, $m\ge 1$, and define a closed subset
of $\mathbb{C}_0$ by
$$
V_{y,m}=\{g\in\mathbb{C}_0:\sup_{t,s\in[0,T]:|t-s|<\frac{T}{m}}|g(t)-g(s)|\ge y).
$$
Next, using the LDP in Remark \ref{R:foll}, we obtain
\begin{align*}
&\limsup_{\varepsilon\rightarrow 0}\varepsilon
\log\mathbb{P}(\sup_{t,s\in[0,T]:|t-s|<\frac{T}{m}}|\widehat{B}_t-\widehat{B}_s|\ge\varepsilon^{-\frac{1}{2}}y)
\le-\inf_{g\in V_{y,m}}J(g) \\
&=-\frac{1}{2}\inf_{g\in V_{y,m}}\inf_{\{f\in\mathbb{H}_0^1:\widehat{f}=g\}}\int_0^T\dot{f}(t)^2dt.
\end{align*}

For a fixed $y> 0$, the sequence 
\begin{equation}
\lambda_m=\inf_{g\in V_{y,m}}\inf_{\{f\in\mathbb{H}_0^1:\widehat{f}=g\}}\int_0^T\dot{f}(t)^2dt,\quad m\ge 1,
\label{E:ku}
\end{equation}
in nondecreasing. It is clear that if we prove that this sequence tends to infinity as $m\rightarrow\infty$, then Corollary
\ref{C:corol} will be established. 

We will next reason by contradiction. Suppose the sequence in (\ref{E:ku}) is bounded. Then there exist sequences $g_m\in V_{y,m}$
and $f_m\in\mathbb{H}_0^1$ with $\widehat{f}_m=g_m$, $m\ge 1$, such that for all $m\ge 1$ and some $C> 0$,
$
\int_0^T\dot{f}_m(t)^2dt\le C.
$
It follows from the compactness statement in Remark \ref{R:remn} that the set $\{g_m\}$ is precompact
in $\mathbb{C}$. By the Arcel\`{a}-Ascoli theorem, the functions $g_m$, $m\ge 1$, are uniformly equicontinuous. 
Recall that for every $m\ge 1$, $g_m\in V_{y,m}$. This means that 
$
\sup_{t,s\in[0,T]:|t-s|<\frac{T}{m}}|g_m(t)-g_m(s)|\ge y,
$
for all $m\ge 1$. It is not hard to see that the previous estimate contradicts the uniform equicontinuity
of the set $\{g_m\}$. It follows that $\lambda_m\rightarrow\infty$ as $m\rightarrow\infty$.

This completes the proof of Corollary \ref{C:corol}.
\subsection{Continuation of the Proof of Theorem \ref{T:27}}\label{SS:conti}
With no loss of generality we can assume that $s_0=1$. 
Using the same ideas as in Section 5 of \cite{G}, we can show that if the term 
$-\frac{1}{2}\varepsilon\int_0^t\sigma(s,\sqrt{\varepsilon}\widehat{B}_s)^2ds$ is removed from the expression for the scaled log-price, then the LDP in Theorem 
\ref{T:27} is not affected. More precisely, this means that it suffices to prove the LDP in Theorem \ref{T:27} for the process 
$\varepsilon\mapsto\widehat{X}^{(\varepsilon)}$, where
\begin{equation}
\widehat{X}^{(\varepsilon)}_t=\int_0^tb(s,\sqrt{\varepsilon}\widehat{B}_s)ds+\sqrt{\varepsilon}
\int_0^t\sigma(s,\sqrt{\varepsilon}\widehat{B}_s)(\bar{\rho}dW_s+\rho dB_s),\quad 0\le t\le T.
\label{E:intl}
\end{equation}

Recall that the functional $\Phi$ associated with the process $\widehat{X}^{(\varepsilon)}$ was defined in (\ref{E:bel}). It will be shown next that the extended contraction principle (see Theorem 4.2.23 in \cite{DZ}) can be applied in our setting. 
Let us first define a sequence of functionals $\Phi_m:\mathbb{C}_0^3\mapsto\mathbb{C}_0$, $m\ge 2$, as follows: 
For $(r,h,l)\in\mathbb{C}_0^3$ and $\frac{jT}{m}<t\le\frac{(j+1)T}{m}$, $1\le j\le m-1$, set
\begin{align}
\Phi_m(r,h,l)(t)&=\frac{T}{m}\sum_{k=0}^{j-1}b\left(\frac{kT}{m},l\left(\frac{kT}{m}\right)\right)+\left(t-\frac{jT}{m}\right)
b\left(\frac{jT}{m},l\left(\frac{jT}{m}\right)\right) \nonumber \\
&\quad+\bar{\rho}\sum_{k=0}^{j-1}
\sigma\left(\frac{kT}{m},l\left(\frac{kT}{m}\right)\right)\left[r\left(\frac{(k+1)T}{m}\right)-r\left(\frac{kT}{m}\right)\right] 
\nonumber \\
&\quad+\bar{\rho}\sigma\left(\frac{jT}{m},l\left(\frac{jT}{m}\right)\right)\left[r\left(t\right)-r\left(\frac{jT}{m}\right)\right] 
\nonumber \\
&\quad+\rho\sum_{k=0}^{j-1}
\sigma\left(\frac{kT}{m},l\left(\frac{kT}{m}\right)\right)\left[h\left(\frac{(k+1)T}{m}\right)-h\left(\frac{kT}{m}\right)\right] 
\nonumber \\
&\quad+\rho\sigma\left(\frac{jT}{m},l\left(\frac{jT}{m}\right)\right)\left[h\left(t\right)-h\left(\frac{jT}{m}\right)\right],
\label{E:lopi}
\end{align}
and for $0\le t\le\frac{T}{m}$, set
$\Phi_m(r,h,l)(t)=tb(0,0)+\bar{\rho}\sigma(0,0)r(t)+\rho\sigma(0,0)h(t)$.
It is not hard to see that for every $m\ge 2$, the mapping $\Phi_m$ is continuous. 

We will next establish that formula (4.2.24) in \cite{DZ} holds in our setting. This formula is used in the 
formulation of the extended contraction principle in \cite{DZ}, Theorem 4.2.23.
\begin{lemma}\label{L:ecp}
For every $c> 0$ and $y> 0$,
$$
\limsup_{m\rightarrow\infty}\sup_{\{(r,f)\in(\mathbb{H}_0^1)^2:\frac{1}{2}\int_0^T\dot{r}(s)^2ds
+\frac{1}{2}\int_0^T\dot{f}(s)^2ds\le c\}}
||\Phi(r,f,\widehat{f})-\Phi_m(r,f,\widehat{f})||_{\mathbb{C}_0}=0.
$$
\end{lemma}

\it Proof. \rm The proof of Lemma \ref{L:ecp} is similar to that of Lemma 21 in \cite{G}. It is not hard to see that
for all $(r,f)\in(\mathbb{H}_0^1)^2$ and $m\ge 2$,
$$
\Phi_m(r,f,\widehat{f})(t)=\int_0^tg_m(s,f)ds+\bar{\rho}\int_0^th_m(s,f)\dot{r}(s)ds+\rho\int_0^th_m(s,f)\dot{f}(s)ds,
$$
where 
$$
g_m(s,f)=\sum_{k=0}^{m-1}b\left(\frac{Tk}{m},\widehat{f}\left(\frac{Tk}{m}\right)\right)
\mathbb{1}_{\{\frac{Tk}{m}\le s\le\frac{T(k+1)}{m}\}},\quad
0\le s\le T,
$$
and
$$
h_m(s,f)=\sum_{k=0}^{m-1}\sigma\left(\frac{Tk}{m},\widehat{f}\left(\frac{Tk}{m}\right)\right)
\mathbb{1}_{\{\frac{Tk}{m}\le s\le\frac{T(k+1)}{m}\}},\quad
0\le s\le T.
$$
Therefore,
\begin{align*}
&\Phi(r,f,\widehat{f})(t)-\Phi_m(r,f,\widehat{f})(t)=\int_0^t[b(s,\widehat{f}(s))-g_m(s,f)]ds \\
&\quad+\bar{\rho}\int_0^t[\sigma(s,\widehat{f}(s))-
h_m(s,f)]\dot{r}(s)ds+\rho\int_0^t[\sigma(s,\widehat{f}(s))-h_m(s,f)]\dot{f}(s)ds.
\end{align*}

For every $\alpha> 0$, denote
$D_{\alpha}=\{w\in\mathbb{H}_0^1:\int_0^T\dot{w}(s)^2ds\le\alpha\}$. It is not hard to see that to prove Lemma \ref{L:ecp}, it suffices to show that for all $\alpha> 0$,
\begin{equation}
\limsup_{m\rightarrow\infty}\left[\sup_{f\in D_{\alpha},w\in D_{\alpha}}\sup_{t\in[0,T]}\left|\int_0^t[\sigma(s,\widehat{f}(s))
-h_m(s,f)]\dot{w}(s)ds\right|\right]=0
\label{E:inter}
\end{equation}
and
\begin{equation}
\limsup_{m\rightarrow\infty}\left[\sup_{f\in D_{\alpha}}\sup_{t\in[0,T]}\left|\int_0^t[b(s,\widehat{f}(s))
-g_m(s,f)]ds\right|\right]=0
\label{E:integl}
\end{equation}

We have
\begin{align}
&\sup_{f\in D_{\alpha},w\in D_{\alpha}}\sup_{t\in[0,T]}\left|\int_0^t[\sigma(s,\widehat{f}(s))
-h_m(s,f)]\dot{w}(s)ds\right| \nonumber \\
&\le\sup_{f\in D_{\alpha},w\in D_{\alpha}}\int_0^T\left|\sigma(s,\widehat{f}(s))
-h_m(s,f)\right||\dot{w}(s)|ds \nonumber \\
&\quad\le\sqrt{T\alpha}\sup_{f\in D_{\alpha}}\sup_{s\in[0,T]}\left|\sigma(s,\widehat{f}(s))
-h_m(s,f)\right|.
\label{E:balk}
\end{align}
We also have
\begin{align}
&\sup_{f\in D_{\alpha}}\sup_{t\in[0,T]}\left|\int_0^t[b(s,\widehat{f}(s))
-g_m(s,f)]ds\right|\le T\sup_{f\in D_{\alpha}}\sup_{s\in[0,T]}\left|b(s,\widehat{f}(s))-g_m(s,f)\right|.
\label{E:rtr}
\end{align}
\begin{lemma}\label{L:compact}
Let $\tau$ be a locally $\omega$-continuous function on $[0,T]\times\mathbb{R}$, $m$ be a positive integer, $f\in\mathbb{H}_0^1$, and
set
$$
q_m(s,f)=\sum_{k=0}^{m-1}\tau\left(\frac{Tk}{m},\widehat{f}\left(\frac{Tk}{m}\right)\right)
\mathbb{1}_{\{\frac{Tk}{m}\le s\le\frac{T(k+1)}{m}\}},\quad
0\le s\le T.
$$
Then,
$
\lim_{m\rightarrow\infty}\sup_{f\in D_{\alpha}}\sup_{s\in[0,T]}\left|\tau(s,\widehat{f}(s))-q_m(s,f)\right|=0.
$
\end{lemma}
 
\it Proof of Lemma \ref{L:compact}. \rm 
The function $\tau$ is locally $\omega$-continuous. We also have the following:
$
\sup_{f\in D_{\alpha}}\sup_{s\in[0,T]}|\widehat{f}(s)|\le M_{\alpha}<\infty.
$
Therefore, for all $f\in D_{\alpha}$ and some $L_{\alpha}> 0$,
\begin{equation}
\sup_{f\in D_{\alpha}}\sup_{s\in[0,T]}|\tau(\widehat{f}(s))-q_m(s,f)|\le L_{\alpha}\omega\left(\frac{T}{m}+
\sup_{f\in D_{\alpha}}E(m,f)\right),
\label{E:fact}
\end{equation}
where 
$
E(m,f)=\sup_{t,u\in[0,T]:|t-u|\le\frac{T}{m}}|\widehat{f}(t)-\widehat{f}(u)|.
$

It follows from the definition of the function $\widehat{f}$ that for all $f\in D_{\alpha}$ and $t,u\in[0,T]$,
$$
|\widehat{f}(t)-\widehat{f}(u)|\le\sqrt{\alpha}\left\{\int_0^T[K(t,v)-K(u,v)]^2dv\right\}^{\frac{1}{2}}
$$
Hence
$$
\sup_{f\in D_{\alpha}}E(m,f)\le\sqrt{\alpha}\sup_{t,u\in[0,T]:|t-u|\le\frac{T}{m}}
\left\{\int_0^T[K(t,v)-K(u,v)]^2dv\right\}^{\frac{1}{2}}\le\sqrt{\alpha}\sqrt{\eta\left(\frac{T}{m}\right)},
$$
for all $m\ge 1$, where $\eta$ is the modulus of continuity in Assumption A. It follows that
$\sup_{f\in D_{\alpha}}E(m,f)\rightarrow 0$ as $m\rightarrow\infty$, and hence (\ref{E:fact}) implies Lemma \ref{L:compact}.

We will next return to the proof of Lemma \ref{L:ecp}. 
It is not hard to see that (\ref{E:inter}) and (\ref{E:integl}) follow from (\ref{E:balk}), (\ref{E:rtr}), and Lemma \ref{L:compact}.

This completes the proof of Lemma \ref{L:ecp}.

It remains to prove that the sequence of processes $\varepsilon\mapsto\Phi_m\left(\sqrt{\varepsilon}W,\sqrt{\varepsilon}B,
\sqrt{\varepsilon}\widehat{B}\right)$ with state space $\mathbb{C}_0$ is an exponentially good approximation to the process $\varepsilon\mapsto\widehat{X}^{(\varepsilon)}$. This property is explained in the next lemma.
\begin{lemma}\label{L:borr}
Suppose that Assumptions A and C hold. Than for every $\delta> 0$,
\begin{align}
&\lim_{m\rightarrow\infty}\,\limsup_{\varepsilon\downarrow 0}\varepsilon\log\mathbb{P}\left(
||\widehat{X}^{(\varepsilon)}-\Phi_m\left(\sqrt{\varepsilon}W,\sqrt{\varepsilon}B,\sqrt{\varepsilon}\widehat{B}\right)||
_{\mathbb{C}_0}>\delta\right)=-\infty.
\label{E:chto}
\end{align}
\end{lemma}

\it Proof of Lemma \ref{L:borr}. \rm It follows from (\ref{E:lopi}) that for all $m\ge 2$,
\begin{align}
&\Phi_m(\sqrt{\varepsilon}W,\sqrt{\varepsilon}B,\sqrt{\varepsilon}\widehat{B})(t)=\int_0^t
b\left(\frac{[msT^{-1}]T}{m},\sqrt{\varepsilon}\widehat{B}_{\frac{[msT^{-1}]T}{m}}\right)ds \nonumber \\
&+\sqrt{\varepsilon}\bar{\rho}\int_0^t\sigma\left(\frac{[msT^{-1}]T}{m},\sqrt{\varepsilon}
\widehat{B}_{\frac{[msT^{-1}]T}{m}}\right)dW_s
+\sqrt{\varepsilon}\rho\int_0^t\sigma\left(\frac{[msT^{-1}]T}{m},\sqrt{\varepsilon}
\widehat{B}_{\frac{[msT^{-1}]T}{m}}\right)dB_s.
\label{E:mysh}
\end{align}
Using (\ref{E:intl}) and (\ref{E:mysh}), we see that in order to prove the equality in (\ref{E:chto}), 
it suffices to show that for every $0<\tau\le 1$,
\begin{equation}
\lim_{m\rightarrow\infty}\,\limsup_{\varepsilon\downarrow 0}\varepsilon\log\mathbb{P}\left(\sqrt{\varepsilon}
\sup_{t\in[0,T]}\left|\int_0^t\sigma_s^{\varepsilon,m}dB_s\right|>\tau\right)=-\infty,
\label{E:tru}
\end{equation}
\begin{equation}
\lim_{m\rightarrow\infty}\,\limsup_{\varepsilon\downarrow 0}\varepsilon\log\mathbb{P}\left(\sqrt{\varepsilon}
\sup_{t\in[0,T]}\left|\int_0^t\sigma_s^{\varepsilon,m}dW_s\right|>\tau\right)=-\infty,
\label{E:trus}
\end{equation}
and
\begin{equation}
\lim_{m\rightarrow\infty}\,\limsup_{\varepsilon\downarrow 0}\varepsilon\log\mathbb{P}\left(
\sup_{t\in[0,T]}\left|\int_0^tb_s^{\varepsilon,m}ds\right|>\tau\right)=-\infty,
\label{E:truss}
\end{equation}
where 
$$
\sigma_s^{\varepsilon,m}=\sigma\left(s,\sqrt{\varepsilon}\widehat{B}_s\right)
-\sigma\left(\frac{[msT^{-1}]T}{m},\sqrt{\varepsilon}\widehat{B}_{\frac{[msT^{-1}]T}{m}}\right), 
\quad 0\le s\le T,\quad m\ge 2,
$$
and
$$
b_s^{\varepsilon,m}=b\left(s,\sqrt{\varepsilon}\widehat{B}_s\right)
-b\left(\frac{[msT^{-1}]T}{m},\sqrt{\varepsilon}\widehat{B}_{\frac{[msT^{-1}]T}{m}}\right), 
\quad 0\le s\le T,\quad m\ge 2.
$$

The formula in (\ref{E:tru}) was obtained in \cite{G1}, formula (53), for $T=1$ and the volatility function of one variable, for which the modulus of continuity $\eta$ in Assumption A is a power function. Lemma 23 in \cite{G1} was used in the proof of formula (53). However, this lemma also holds in the case of a general modulus of continuity satisfying Assumption A. This was established in Corollary \ref{C:corol} in the present paper. The proof of formula (53) in \cite{G1} is rather complicated, but it is not very difficult to adapt it to our setting thus establishing the equalities in (\ref{E:tru}) and (\ref{E:trus}). We will sketch below how this adaptation can be performed, and leave filling in the details as an exercise for the interested reader. First, we will prove the equality in (\ref{E:truss}), using the same techniques as in the proof of (53) in \cite{G1}, and by taking into account 
Corollary \ref{C:corol} in the present paper instead of Lemma 23 in \cite{G1}. The proof of (\ref{E:truss}) is simpler than that of 
(\ref{E:tru}) and (\ref{E:trus}), since no estimates of stochastic integrals are needed in it. 
We include the proof of (\ref{E:truss}) below to illustrate the main ideas used in the proof of (\ref{E:tru}) and (\ref{E:trus}).

Recall that the class $MC$ was defined in Remark \ref{R:let}. The next simple statement will be used in the proof of (\ref{E:truss}).
\begin{lemma}\label{L:auxi}
For every $\mu\in MC$ and $N\in MC$, there exists a positive strictly decreasing continuous function $q(r)$, $0< r< r_0$, such that 
the following conditions hold: $\lim_{r\rightarrow 0}q(r)=\infty$ and $\lim_{r\rightarrow 0}N\left(q(r)\right)\mu\left(r\right)=0$.
\end{lemma}

\it Proof. \rm Define a function by $q(r)=N^{-1}\left(\mu(r)^{-\frac{1}{2}}\right)$, $0< r< r_0$. It is easy to see that the function $q$ satisfies the conditions in Lemma \ref{L:auxi}.

We will next return to the proof of (\ref{E:truss}). Recall that 
we assumed that $\omega$ and $L$ are functions from the class $MC$ (see Remark \ref{R:let}).

Let $q$ be the function obtained by applying Lemma \ref{L:auxi} to the moduli of continuity $\mu(r)=\sqrt{\omega(2r)}$ 
and $N(r)=L(2r)$. Hence, $q$ is a positive strictly decreasing continuous function on $(0,r_0)$ such that
$\displaystyle{\lim_{r\rightarrow 0}q(r)=\infty}$ and
\begin{equation}
\lim_{r\rightarrow 0}L\left(2q(r)\right)\sqrt{\omega\left(2r\right)}=0.
\label{E:tyty}
\end{equation}

For every $m\ge 2$, $0<\varepsilon< 1$, and $0<r< r_0$, define a random variable by 
\begin{equation}
\xi^{\varepsilon,m,r}(w)=\inf_{s\in[0,T]}\left\{\sqrt{\varepsilon}\left[\frac{r}{q(r)}|\widehat{B}_s(w)|+
|\widehat{B}_s(w)-\widehat{B}_{\frac{[msT^{-1}]T}{m}}(w)|\right]> r\right\}.
\label{E:tt1}
\end{equation}
If for some $w\in\Omega$, the set appearing on the right-hand side of (\ref{E:tt1}) is empty, 
we put $\xi^{\varepsilon,m,r}(w)=T$.
It is not hard to see that $\xi^{\varepsilon,m,r}$ is an ${\cal F}^B_t$-stopping time. Here we use the fact that 
the filtration ${\cal F}^B_t$ is right-continuous and the process $\widehat{B}$ is continuous. 
We also have
\begin{equation}
\sqrt{\varepsilon}\max\left[|\widehat{B}_s|,
|\widehat{B}_{\frac{[msT^{-1}]T}{m}}|\right]\le q(r)\quad\mbox{and}\quad\sqrt{\varepsilon}
|\widehat{B}_s-\widehat{B}_{\frac{[msT^{-1}]T}{m}}|\le r,
\label{E:tt2}
\end{equation}
for all $s\le\xi^{\varepsilon,m,r}$. Moreover, it follows from (\ref{E:tt2}) and the $\omega$-continuity of the function $\sigma$ that
\begin{equation}
|\sigma_s^{\varepsilon,m}|\le L\left(T+q(r)\right)\omega\left(\frac{T}{m}+r\right),
\label{E:urah}
\end{equation}
if $s\le\xi^{\varepsilon,m,r}$. It is clear that for any $r<r_0$, there exists a positive integer $m_r$ such that 
\begin{equation}
|\sigma_s^{\varepsilon,m}|\le L\left(2q(r)\right)\omega\left(2r\right),\quad m>m_r.
\label{E:tr0}
\end{equation}
It is also clear from (\ref{E:tyty}) 
\begin{equation}
\lim_{r\rightarrow 0}L\left(2q(r)\right)\omega\left(2r\right)=0.
\label{E:lpl}
\end{equation}
Moreover, for every $r> 0$,
\begin{equation}
\mathbb{P}\left(\sup_{t\in[0,T]}\left|\int_0^tb_s^{\varepsilon,m}ds\right|>\tau\right)
\le\mathbb{P}\left(\sup_{t\in[0,\xi^{\varepsilon,m,r}]}\left|\int_0^tb_s^{\varepsilon,m}ds\right|>\tau\right)
+\mathbb{P}(\xi^{\varepsilon,m,r}< T).
\label{E:tt3}
\end{equation}
It follows from (\ref{E:tt2}) and the $\omega$-continuity condition for $b$ that 
\begin{align*}
&\mathbb{P}\left(\sup_{t\in[0,\xi^{\varepsilon,m,r}]}\left|\int_0^tb_s^{\varepsilon,m}ds\right|>\tau\right) \\
&\le\mathbb{P}\left(TL\left(T+q(r)\right)\sqrt{\omega\left(\frac{T}{m}+r\right)}
\sup_{t\in[0,\xi^{\varepsilon,m,r}]}\sqrt{\omega\left(\frac{T}{m}+\sqrt{\varepsilon}
|\widehat{B}_t-\widehat{B}_{\frac{[mtT^{-1}]T}{m}}|\right)}>\tau\right) \\
&\le\mathbb{P}\left(TL\left(T+q(r)\right)\sqrt{\omega\left(\frac{T}{m}+r\right)}
\sqrt{\omega\left(\frac{T}{m}+\sqrt{\varepsilon}
\sup_{t\in[0,T]}|\widehat{B}_t-\widehat{B}_{\frac{[mtT^{-1}]T}{m}}|\right)}>\tau\right)  \\
&\le\mathbb{P}\left(\sqrt{\varepsilon}
\sup_{t\in[0,T]}|\widehat{B}_t-\widehat{B}_{\frac{[mtT^{-1}]T}{m}}|
>\omega^{-1}\left(\tau^2T^{-2}L\left(T+q(r)\right)^{-2}\omega\left(\frac{T}{m}+r\right)^{-1}\right)
-\frac{T}{m}\right).
\end{align*}
Since (\ref{E:tyty}) holds, we can find $r_0> 0$ such that for all $0< r< r_0$ and $m\ge 1$,
$$
\omega^{-1}\left(\tau^2T^{-2}L\left(2q(r)\right)^{-2}\omega\left(2r\right)^{-1}\right)
-\frac{T}{m}> y_0,
$$ 
where $y_0$ is a number appearing in Corollary \ref{C:corol}. Fix $r$ with $0< r< r_0$. Then there exists a positive integer $m_r$
such that for all $m\ge m_r$,
\begin{align*}
&\omega^{-1}\left(\tau^2T^{-2}L\left(T+q(r)\right)^{-2}\omega\left(\frac{T}{m}+r\right)^{-1}\right) 
-\frac{T}{m} \\
&\ge\omega^{-1}\left(\tau^2T^{-2}L\left(2q(r)\right)^{-2}\omega\left(2r\right)^{-1}\right)
-\frac{T}{m}> y_0
\end{align*}
Next, using (\ref{E:cvf}) in Corollary \ref{C:corol}, we see that for every $0< r< r_0$,
\begin{align}
&\limsup_{m\rightarrow\infty}\limsup_{\varepsilon\rightarrow 0}\varepsilon
\log\mathbb{P}\left(\sup_{t\in[0,\xi^{\varepsilon,m,r}]}\left|\int_0^tb_s^{\varepsilon,m}ds\right|>\tau\right)=-\infty.
\label{E:megin}
\end{align}

We will next estimate the second term on the right-hand side of (\ref{E:tt3}). It is not hard to see that the following inequality holds:
\begin{align}
&\mathbb{P}(\xi^{\varepsilon,m,r}< T) \nonumber \\
&\le\mathbb{P}\left(\sqrt{\varepsilon}\sup_{t\in[0,T]}|\widehat{B}_t|>\frac{q(r)}{2}\right)
+\mathbb{P}\left(\sqrt{\varepsilon}\sup_{t\in[0,T]}|\widehat{B}_s-\widehat{B}_{\frac{[msT^{-1}]T}{m}}|>\frac{r}{2}\right).
\label{E:ptu1}
\end{align}
Using the exponential estimate for the distribution function of the supremum of a Gaussian process (see the reference in the proof of 
(36) in \cite{G1}), we get
$$
\mathbb{P}\left(\sqrt{\varepsilon}\sup_{t\in[0,T]}|\widehat{B}_t|>\frac{q(r)}{2}\right)\le\exp\left\{-C_1\frac{q(r)^2}{\varepsilon}\right\},\,\,r< r_1,\,\,\varepsilon<\varepsilon_0,
$$
for some constant $C_1> 0$. Therefore, since $q(r)\rightarrow\infty$ as $r\rightarrow 0$,
\begin{equation}
\limsup_{r\rightarrow 0}\limsup_{\varepsilon\rightarrow 0}\varepsilon
\log\mathbb{P}\left(\sqrt{\varepsilon}\sup_{t\in[0,T]}|\widehat{B}_t|>\frac{q(r)}{2}\right)=-\infty.
\label{E:ryt}
\end{equation}
In addition, using Corollary \ref{C:corol}, we see that for all $r> 0$,
\begin{equation}
\limsup_{m\rightarrow\infty}\limsup_{\varepsilon\rightarrow 0}\varepsilon\log
\mathbb{P}\left(\sqrt{\varepsilon}\sup_{t\in[0,T]}|\widehat{B}_s-\widehat{B}_{\frac{[msT^{-1}]T}{m}}|>\frac{r}{2}\right)=-\infty.
\label{E:ptu3}
\end{equation}

The following inequality holds for all $a> 0$ and $b> 0$,
\begin{equation}
\log(a+b)\le\max\{\log(2a),\log(2b)\}.
\label{E:ppi}
\end{equation}
Next, using (\ref{E:ppi}), we can glue together the estimates in (\ref{E:ryt}) and (\ref{E:ptu3}). This allows us to obtain from (\ref{E:ptu1})
that
\begin{align}
\limsup_{r\rightarrow 0}\limsup_{m\rightarrow\infty}\limsup_{\varepsilon\rightarrow 0}
\varepsilon\log\mathbb{P}(\xi^{\varepsilon,m,r}< T)=-\infty.
\label{E:ali1}
\end{align}
 
Finally, by taking into account (\ref{E:tt3}), (\ref{E:megin}), (\ref{E:ali1}), and (\ref{E:ppi}), we establish (\ref{E:truss}).  

We will next prove (\ref{E:tru}). The proof of (\ref{E:trus}) is similar, but simpler, than that of (\ref{E:tru}) since the processes 
$\widehat{B}$ and $W$ are independent. Observe that 
\begin{equation}
\mathbb{P}\left(\sup_{t\in[0,T]}\left|\int_0^t\sigma_s^{\varepsilon,m}dB_s\right|>\tau\right)
\le\mathbb{P}\left(\sup_{t\in[0,\xi^{\varepsilon,m,r}]}\left|\int_0^t\sigma_s^{\varepsilon,m}dB_s\right|>\tau\right)
+\mathbb{P}(\xi^{\varepsilon,m,r}< T).
\label{E:tt33}
\end{equation}
We already know from the proof of (\ref{E:truss}) how to handle the second term on the right-hand side of (\ref{E:tt33}). It remains to estimate the first term.

It will be shown next that the process $N^{\varepsilon,m,r}(t)=\int_0^{t\wedge\xi^{\varepsilon,m,r}}\sigma_s^{\varepsilon,m}dB_s$, 
$t\in[0,T]$, is an ${\cal F}^B_t$-martingale. Although the proof of the previous statement is standard, we decided to include it below for the sake of convenience. It is not hard to see that the process $M(t)=\int_0^t\sigma_s^{\varepsilon,m}dB_s$, $t\in[0,T]$, is a local martingale. Let $\tau_n\uparrow T$ be a localizing sequence of stopping times for the process $M$. Since for every $n\ge 1$, the process $M_n(t)=
M(t\wedge \tau_n)$, $t\in[0,T]$, is a martingale, the process $M_n^{(\varepsilon,m,r)}(t)=M(t\wedge\tau_n\wedge\xi^{\varepsilon,m,r})$, $t\in[0,T]$, is also a martingale (see Corollary 3.6 in \cite{RY}). Hence for all $0\le s\le t\le T$,
\begin{equation}
\mathbb{E}[M_n^{(\varepsilon,m,r)}(t)|{\cal F}_s^B]=M_n^{(\varepsilon,m,r)}(s)
\label{E:lp}
\end{equation}
$\mathbb{P}$-a.s. on $\Omega$. Our next goal is to pass to the limit as $n\rightarrow\infty$ in (\ref{E:lp}). 
The process $M$ is a continuous  process. Therefore for every $u\in[0,T]$, 
\begin{equation}
\lim_{n\rightarrow\infty}M_n^{(\varepsilon,m,r)}(u)\rightarrow N^{\varepsilon,m,r}(u)
\label{E:sss}
\end{equation} 
$\mathbb{P}$-a.s. on $\Omega$. To justify the possibility of passing to the limit under the conditional expectation sign on the left-hand side of (\ref{E:lp}), we use the inequality
\begin{equation}
\mathbb{E}[\sup_{n\ge 1}|M_n^{(\varepsilon,m,r)}(t)|]<\infty,\,\,t\in[0,T],
\label{E:pl}
\end{equation}
and the dominated convergence theorem. It remains to prove (\ref{E:pl}). It follows from Doob's maximal inequality and (\ref{E:urah}) that for every $n\ge 1$,
\begin{align*} 
&\mathbb{E}[\sup_{0\le u\le T}M_n^{(\varepsilon,m,r)}(u)^2]\le 4\mathbb{E}[M_n^{(\varepsilon,m,r)}(T)^2]
\le 4\int_0^{\xi^{\varepsilon,m,r}}\left(\sigma_s^{\varepsilon,m}\right)^2ds
\\
&\le 4T L\left(T+q(r)\right)^2\omega\left(\frac{T}{m}+r\right)^2.
\end{align*}
It is not hard to show that $n\mapsto\sup_{0\le u\le T}M_n^{(\varepsilon,m,r)}(u)^2$ is an increasing sequence of random variables.
Moreover,
$
\sup_{0\le u\le T}N^{(\varepsilon,m,r)}(u)^2\le\lim_{n\rightarrow\infty}\sup_{0\le u\le T}M_n^{(\varepsilon,m,r)}(u)^2.
$
It follows from the monotone convergence theorem that
\begin{equation}
\mathbb{E}[\sup_{0\le u\le T}N^{(\varepsilon,m,r)}(u)^2]\le 4T L\left(T+q(r)\right)^2\omega\left(\frac{T}{m}+r\right)^2<\infty.
\label{E:opi}
\end{equation}
Since for all $t\in[0,T]$,  
$\sup_{n\ge 1}|M_n^{(\varepsilon,m,r)}(t)|\le\sup_{0\le u\le T}|N^{(\varepsilon,m,r)}(u)|$, and (\ref{E:opi}) holds, we get 
(\ref{E:pl}). Finally, we can show that the process $N^{(\varepsilon,m,r)}$ is a martingale, by using (\ref{E:lp}), (\ref{E:sss}),
(\ref{E:pl}).

Let us continue the proof of (\ref{E:tru}). Reasoning as in the proof of formulas (61) and (62) in \cite{G1}, and taking into account 
(\ref{E:tr0}), we see that for every $\lambda> 0$, the stochastic exponential
$$ 
{\cal E}_t=\exp\left\{-\frac{1}{2}\lambda^2\varepsilon\int_0^{t\wedge\xi^{\varepsilon,m,r}}(\sigma_s^{\varepsilon,m})^2ds
+\lambda\sqrt{\varepsilon}\int_0^{t\wedge\xi^{\varepsilon,m,r}}\sigma_s^{\varepsilon,m}dB_s\right\}
$$
is a martingale. Moreover, for all $\lambda> 0$, $r< r_0$, and $m\ge m_r$, we have
$$
\mathbb{E}\left[\exp\left\{\lambda\sqrt{\varepsilon}\int_0^{\xi^{\varepsilon,m,r}}\sigma_s^{\varepsilon,m}dB_s\right\}\right]
\le\exp\left\{\frac{T}{2}\lambda^2\varepsilon L(q(r))^2\omega(r)^2\right\}
$$
(see the proof of formulas (61) and (62) in \cite{G1}). Hence the process
$$
t\mapsto\exp\left\{\lambda\sqrt{\varepsilon}\int_0^{t\wedge\xi^{\varepsilon,m,r}}\sigma_s^{\varepsilon,m}dB_s\right\}
=\exp\left\{\lambda\sqrt{\varepsilon}\int_0^t\mathbb{1}_{s\le\xi^{\varepsilon,m,r}}\sigma_s^{\varepsilon,m}dB_s\right\}
$$
is a positive submartingale. Here we refer to Proposition 3.6 in Chapter 1 of \cite{KaS}. Therefore, we can use the first submartingale inequality in \cite{KaS} (see Theorem 3.8 in Chapter 1 of \cite{KaS}). This gives
\begin{align*}
&\mathbb{P}\left(\sup_{t\in[0,T]}\exp\left\{\lambda\sqrt{\varepsilon}
\int_0^t\mathbb{1}_{s\le\xi^{\varepsilon,m,r}}\sigma_s^{\varepsilon,m}dB_s\right\}\ge e^{\lambda\tau}\right) \\
&\le\mathbb{E}\left[\exp\left\{\lambda\sqrt{\varepsilon}\int_0^{\xi^{\varepsilon,m,r}}
\sigma_s^{\varepsilon,m}dB_s-\lambda\tau\right\}\right]\le\exp\left\{\frac{T}{2}\lambda^2\varepsilon L(q(r))^2\omega(r)^2
-\lambda\tau\right\},
\end{align*}
and it follows that
\begin{align*}
&\mathbb{P}\left(\sup_{t\in[0,\xi^{\varepsilon,m,r}]}\sqrt{\varepsilon}
\int_0^t\sigma_s^{\varepsilon,m}dB_s\ge\tau\right)
\le\exp\left\{\frac{T}{2}\lambda^2\varepsilon L(q(r))^2\omega(r)^2-\lambda\tau\right\}.
\end{align*}
Similarly, we prove that
\begin{align*}
&\mathbb{P}\left(\sup_{t\in[0,\xi^{\varepsilon,m,r}]}\sqrt{\varepsilon}(-1)
\int_0^t\sigma_s^{\varepsilon,m}dB_s\ge\tau\right)
\le\exp\left\{\frac{T}{2}\lambda^2\varepsilon L(q(r))^2\omega(r)^2-\lambda\tau\right\},
\end{align*}
and hence
\begin{align*}
&\mathbb{P}\left(\sup_{t\in[0,\xi^{\varepsilon,m,r}]}\sqrt{\varepsilon}
\left|\int_0^t\sigma_s^{\varepsilon,m}dB_s\right|\ge\tau\right)
\le 2\exp\left\{\frac{T}{2}\lambda^2\varepsilon L(q(r))^2\omega(r)^2-\lambda\tau\right\}.
\end{align*}
Next, using the previous inequality and (\ref{E:lpl}), we obtain
\begin{equation}
\limsup_{r\rightarrow 0}\limsup_{m\rightarrow\infty}\limsup_{\varepsilon\rightarrow 0}\varepsilon\log\mathbb{P}\left(\sup_{t\in[0,\xi^{\varepsilon,m,r}]}\sqrt{\varepsilon}
\left|\int_0^t\sigma_s^{\varepsilon,m}dB_s\right|\ge\tau\right)=-\infty.
\label{E:erw1}
\end{equation}
Finally, we can finish the proof of (\ref{E:tru}) using (\ref{E:tt33}) and (\ref{E:erw1}) and reasoning as in the second half of the proof
of (\ref{E:truss}).

The proof of Lemma \ref{L:borr} is thus completed.

We will next return to the proof of Theorem \ref{T:27}. By taking into account Theorem \ref{T:ladp}, Lemmas \ref{L:ecp} and \ref{L:borr}, and applying the extended contraction 
principle (see Theorem 4.2.23 in \cite{DZ}), we show that the process $\varepsilon\mapsto\widehat{X}^{(\varepsilon)}$ satisfies the large deviation principle 
with speed $\varepsilon^{-1}$ and good rate function $\widehat{Q}_T$ (see the definition in (\ref{E:vunzi})). 

This completes the proof of Theorem \ref{T:27}.

\section{Acknowledgments}\label{S:ack}
I thank Frederi Viens for valuable suggestions. I also gratefully acknowledge
the detailed comments of the anonymous referees.

\end{document}